\documentclass[10pt,onecolumn,superscriptaddress,showpacs,aps,pre]{revtex4-1}
\usepackage{color}
\usepackage{amsmath}
\usepackage{CJK}
\usepackage{psfrag}
\usepackage{graphicx}
\usepackage{subfig}

\begin{document}

\title{The dynamics of ensemble of neuron-like elements with excitatory couplings}
\author{Alexander G. Korotkov}
\affiliation{Control Theory Department, Institute of Information Technologies, Mathematics and Mechanics, Lobachevsky University, Gagarin ave. 23, Nizhny Novgorod, 603950, Russia}
\author{Alexey O. Kazakov}
\affiliation{Faculty of Informatics, Mathematics, and Computer Science, National Research University Higher School of Economics, Bolshaja Pecherskaja Str. 5/12, Nizhny Novgorod, 603155, Russia}
\author{Tatiana A. Levanova}
\affiliation{Control Theory Department, Institute of Information Technologies, Mathematics and Mechanics, Lobachevsky University, Gagarin ave. 23, Nizhny Novgorod, 603950, Russia}
\author{Grigory V. Osipov}
\affiliation{Control Theory Department, Institute of Information Technologies, Mathematics and Mechanics, Lobachevsky University, Gagarin ave. 23, Nizhny Novgorod, 603950, Russia}
\date{\today}

\begin{abstract}
We study the phenomenological model of ensemble of two FitzHugh-Nagumo neuron-like elements with symmetric excitatory couplings. The main advantage of proposed model is the new approach to model of coupling which is implemented by smooth function that approximate rectangular function. The proposed coupling depends on three parameters that define the beginning of activation of an element $\alpha$, the duration of the activation $\delta$ and the strength of the coupling $g$. We observed a rich diversity of types of neuron-like activity, including regular in-phase, anti-phase and sequential spiking activities. In the phase space of the system, these regular regimes correspond to specific asymptotically stable periodic motions (limit cycles). We also observed a  chaotic anti-phase activity, which corresponds to a strange attractor that appears due to the cascade of period doubling bifurcations of limit cycles.

We also provide the detailed study of bifurcations which lead to transitions between all these regimes and detect on the $(\alpha, \delta)$ parameter plane those regions that correspond to  the above-mentioned regimes. We also show numerically the existence of bistability regions when various non-trivial regimes coexist. For example, in some regions, one can observe either anti-phase or in-phase oscillations depending on initial conditions. We also specify regions corresponding to coexisting various types of sequential activity.
\end{abstract}

\maketitle

\section{Introduction}
Understanding the processes of brain and the nervous system functioning is one of the main problems in contemporary science \cite{Markram}. The brain contains networks connecting approximately $10^{12}$ neurons with $10^{15}$ synapses, which are responsible for sensor and motor functions, cognitive activity, emotions and, as a result, they define behavior and consciousness of individual \cite{DAngelo}. Neurons are excitable cells that can generate action potentials. The couplings between neurons are organized using special contacts called synapses \cite{Synapses}. Depending on the type of impact, it is possible to divide synapses into excitatory and inhibitory ones. In a similar way, depending on the method of signal transmission all synapses can be divided into electrical, in which the signals are transmitted by electric current, and chemical, in which signal transmission is carried out by means of a physiologically active substance called the neurotransmitter \cite{Neuroscience}.

Despite numerous investigations and significant progress in neuroscience, the problem of studying the features of neural connections and their influence on the functioning of the brain is still far from its solution. To solve it, one needs not only experimental studies, but also the development of various mathematical methods and models, beginning with the cellular level and ending with the modelling of the brain as a single large system. At the present moment, one can find a number of various models that allow to study molecular and cellular processes in neurons, as well as the dynamics of individual neurons and neural ensembles \cite{Models1}. In the most general context, all mathematical models can be divided into two classes: realistic biological models and phenomenological modifications \cite{Models1}. In the first case, the various biological data and the biophysical principles \cite{Details} should be taken into account first of all. For example, the most famous such models are the Hodgkin-Huxley model \cite{HH} as well as its various modifications \cite{Shilnikov, Shilnikov1, Shilnikov2}. In the second case, only reproduction of a specific biological phenomenon observed in the experiment is required, and the model is constructed without consideration of certain biological details. Some phenomenological models are reductions of the original realistic biological models, for example, as well-known FitzHugh-Nagumo model \cite{FHN}. In the present paper we will consider the model of neuron ensemble of such a type.

In this case, for modeling, it can be assumed that coupling due to which the action potential (the current $I_{syn}(t)$) is transmitted from the presinaptic neuron to the postsynaptic one is determined by the following equation:
\begin{equation*}
I_{syn}(t) = -g(t)\cdot(V(t) - V_{rev}).
\end{equation*}
Here $V_{rev}$ is a reversal potential, $g(t)$ is a synaptic conductance, $V(t)$ is a membrane potential of a postsynaptic neuron \cite{Roth-Rossum, Baladron-Fasoli-Faugeras-Touboul}. The transmission of action potential can be mathematically modelled in different ways, for example, using differential equation \cite{Gao}, equations with delayed argument \cite{Mobius, Rankovic} or the so-called $\alpha$-function \cite{Shin-Lee-Kim, Hansel-Mato-Meuner}. In the last approach, $I_{syn}(t)$ is calculated from the following relation:
\begin{equation*}
I_{syn}(t) = -g\cdot\alpha(t - t_0)\cdot(V(t) - V_{rev}),
\end{equation*}
where $g$ is peak synaptic conductance, $t_0$ is a transmission start time, $\alpha(t) = \frac{t}{\tau}\exp(-\frac{t}{\tau})$ is $\alpha$-function with parameter $\tau$ that define characteristic time of interaction between pre- and postsynaptic neurons.

It is worth noting that the FitzHugh-Nagumo model is widely used as the basic model of the individual element in simulating of neural ensembles, thus, further, we give a short review on this topic, paying more attention to the modelling of minimal ensembles. First of all we note the paper \cite{Doss-Bachelet-Francoise-Piquet}, where the occurrence of bursting activity in such ensembles was investigated. Recall that under bursting activity here is meant periodic or chaotic oscillations for which sequences of fast oscillations (spikes) alternate with slow subthreshold motions \cite{Izhikevich2000}. The papers \cite{Komarov-Osipov-Suykens-Rabinovich, Rabinovich} are devoted to the investigation of sequential activity --- such a type of neuron-like regimes for which the activity from one (active) element is transmitted to another element, while active element become suppressed for some time. Also it is important to note some papers in which synchronization phenomena in minimal ensembles were investigated, see e.g. \cite{Nguyen, Shin-Lee-Kim}. Concerning ensembles of FitzHugh-Nagumo elements with more complex coupling, it is important to note here the papers \cite{Binczak, Jacquir}, where various non-trivial (and even chaotic) neuron-like regimes were studied for such ensembles with unidirectional couplings. Similar results were also obtained for the bi-directional coupling in \cite{Hoff}. We also note the papers \cite{Campbell, Tehrani, Yanagita} where it was shown that the dynamics in minimal ensembles with synaptic couplings can exhibit multistability -- when  different attractors coexist in the phase space of the ensemble and type of the stable dynamics in the system depends on the initial conditions. It is also important to refer the papers devoted to the investigation of influence of different factors, such as noise perturbation \cite{Brown-Feng-Feerick} and perturbation by external stimulus \cite{Wang-Zhang-Deng}, as well as influence of time delay \cite{Song, Nguyen} and topology of the couplings \cite{Wang-Tian-Dhamala-Liu,Hoff} on the dynamics of two coupled FitzHugh-Nagumo systems.

The main goal of this paper is to study both various types of neuron-like activity, such as the in-phase and anti-phase spiking regimes, various regimes of sequential activity, and also bifurcation scenarios of the appearance and destruction of these regimes in the model of two excitatory coupled FitzHugh-Nagumo elements. Let us explain what is meant here by the listed regimes of neuron-like activity. Here the \textit{in-phase regime} is a regime in which the state of ensemble changes periodically, and the states of both elements coincide. Time series for this regime are shown in Fig. \ref{time_serieses_b}. The \textit{anti-phase regime} is a regime in which the state of the first element approximately coincides with the state of the second one shifted by half period: $x_1(t) \approx x_2 (t\pm T/2), y_1(t) \approx y_2 (t\pm T/2)$, where $T$ is the period of the limit cycle. Time series for this regime are shown in Fig. \ref{time_serieses_a}. The \textit{regime of sequential activity} is a regime at which the activity is switching between elements, so that the elements sequentially, without a delay, generate action potentials and then remain suppressed for a while (examples of various regimes of sequential activity are given in the Figs. \ref{time_serieses_c} , \ref{time_serieses_d}, \ref{time_serieses_e}, \ref{time_serieses_f}).

In this paper, we propose a new phenomenological model of the synaptic coupling between elements that is simple from the computational point of view. This coupling, nevertheless, allows to obtain a wide variety of regimes of neuron-like activity observed in experiments and biologically plausible models. It is important to note that in our modelling we use a coupling in which an impulse from the active presynaptic element arrives at postsynaptic element when the polar angle of the presynaptic element lies in the range specified by the two governing parameters of the model. This allow us to model some important properties of the synapse dynamics.

\section{The model}
The FitzHugh-Nagumo system is one of the simplest models that describe the behavior of a neuron-like element. It can be written in the following form
\begin{equation} \label{1_element}
\begin{cases}
\epsilon \mathop{x}\limits^\cdot = x - x^3/3 - y\\
\mathop{y}\limits^\cdot = x - a
\end{cases}
\end{equation}
Here $x$ is a variable describing the dynamics of membrane potential of the neuron-like element, and $y$ is a recovery variable, also $\epsilon$ is  small parameter, $0 < \epsilon << 1$. Thus, system \eqref{1_element} is a typical slow-fast system with the fast variable $x$ and slow variable $y$.

\begin{figure}[h!]
	\center{
	\subfloat[Phase portait, $-1 < a < 1$]
	{
		\includegraphics[width = 0.3\linewidth]{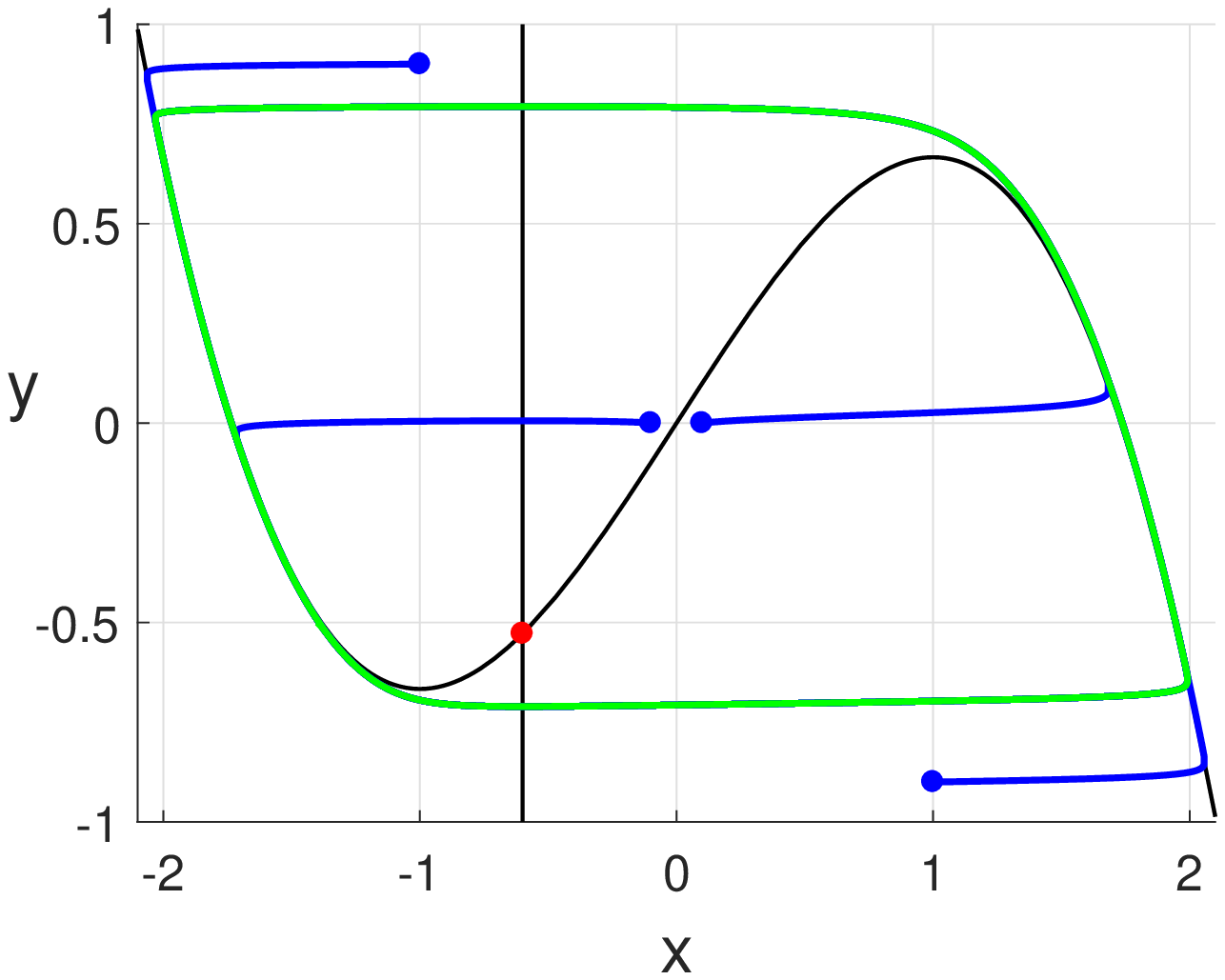}
		\label{1_el_a}
	}
	\hspace{1cm}
	\subfloat[Time series for $x(t)$, $-1 < a < 1$]
	{
		\includegraphics[width = 0.3\linewidth]{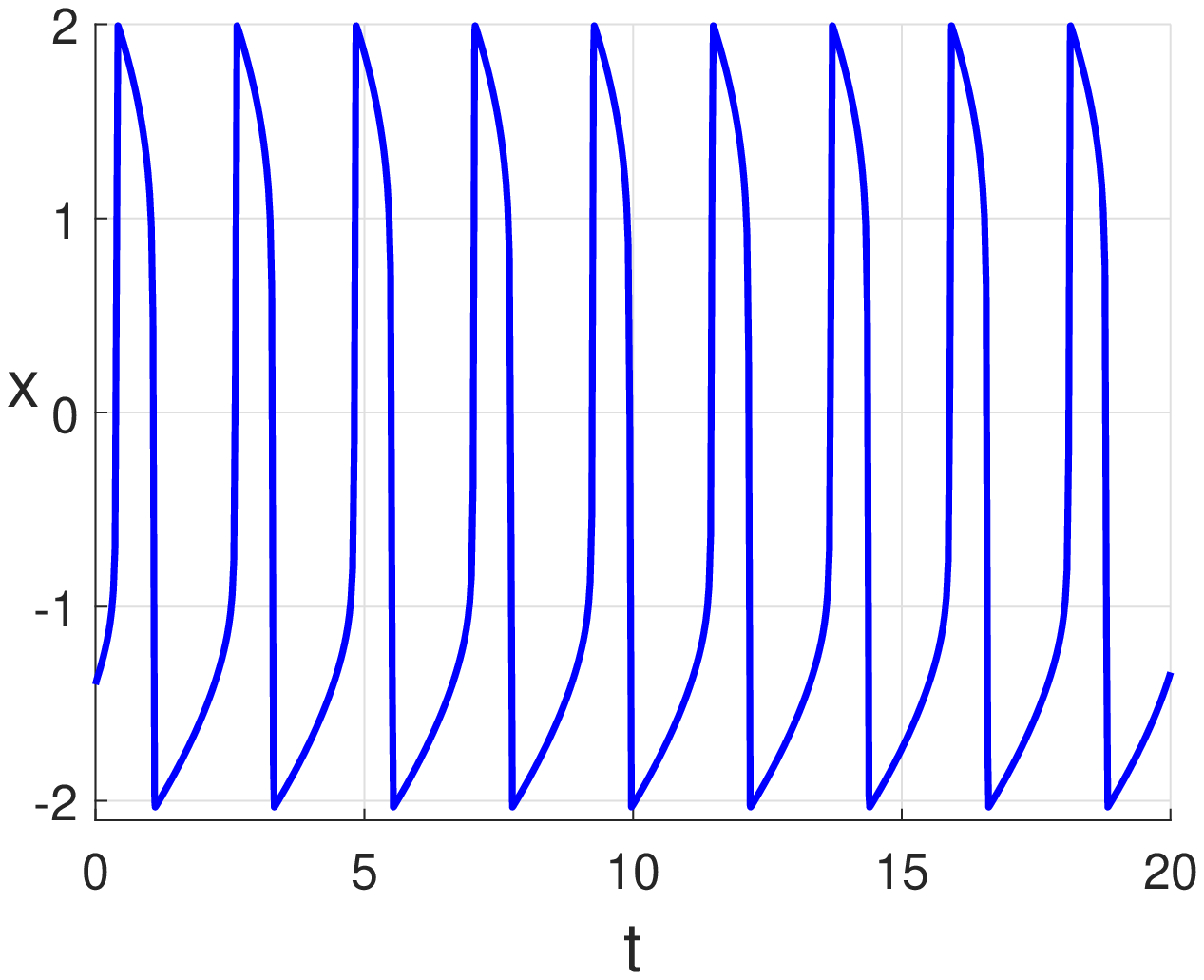}
		\label{1_el_b}
	}\\
	\subfloat[Phase portrait, $a < -1$]
	{
		\includegraphics[width = 0.3\linewidth]{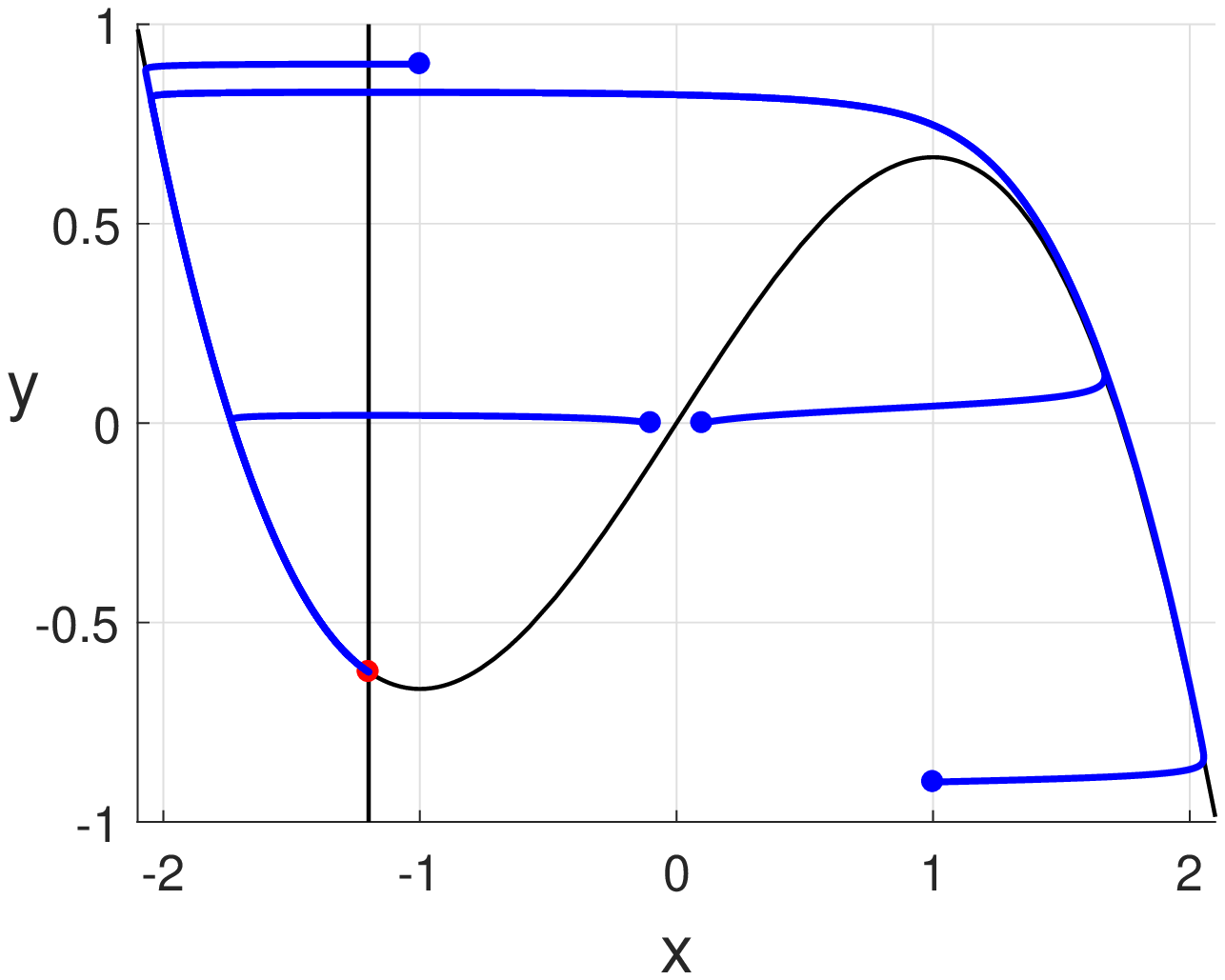}
		\label{1_el_c}
	}
	\hspace{1cm}
	\subfloat[External impulse (bottom) and the responce for it (top) for $a < -1$]
	{
		\includegraphics[width = 0.3\linewidth]{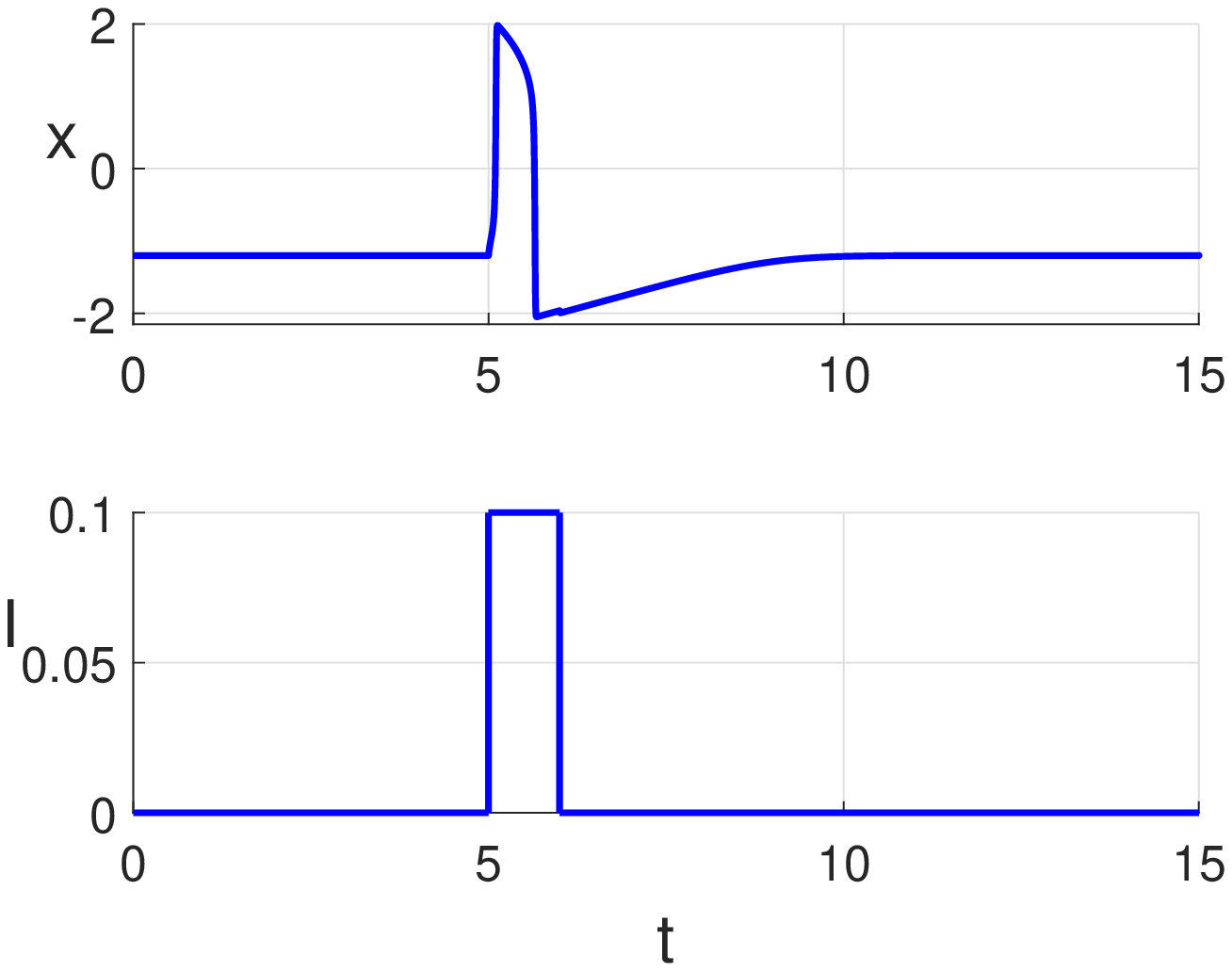}
		\label{1_el_d}
	}
}
	\caption{(a) Phase portrait of system \eqref{1_element} in the self-oscillatory regime. (b) Time series for the self-oscillatory regime. (c) Phase portrait of the system \eqref{1_element} in the excitable regime. (d) The response of excitable element to external impulse. In (a) and (c) isoclines $y = x - x^3/3$ and $x = a$ are shown in black colour, stable limit cycle is shown in green colour, trajectories of the system \eqref{1_element} with its initial points are shown in blue colour; the red point is the equilibrium.}
	\label{1_el}
\end{figure}

For any values of the parameters $a$ and $\epsilon$, system \eqref{1_element} has only one equilibrium $\tilde{O}$ $(a, a - a^3/3)$. It is asymptotically stable for $\left | a \right | \geq 1$ and completely unstable for $-1 < a < 1$. At $a = \pm 1$ this equilibrium undergoes the supercritical Andronov-Hopf bifurcation. As a result, for $|a| > 1$ the stable limit cycle appears in the phase space of the system.

In the Fig. \ref{1_el} we show examples of phase portraits of the system \eqref{1_element} in the case $ a = -0.6 $ (Fig. \ref{1_el_a}) and in the case of $ a = -1.2 $ (Fig. \ref{1_el_c}). The lines corresponding to the isoclines $x-x^3/3-y=0 $ and $x-a=0$ are shown in black colour. In Fig. \ref{1_el_a} green colour marks a stable limit cycle, which is a global attractor. Corresponding time series are shown in Fig. \ref{1_el_b}. As one can see from Fig. \ref{1_el_b}, the oscillations include fast and slow changes of the state of system. In the phase portrait, the slow motions corresponds to the passage of the phase point near the isocline $x-x^3/3-y=0$, since in a sufficiently small neighbourhood of this curve the values of the fast variable $x$ are very close to a constant. The fast motions correspond to motions along trajectories close to the horizontal lines. In this case the vertical component of the phase velocity vector $\mathop{y}\limits^\cdot$ is negligible compared to the horizontal component $\mathop {x}\limits^\cdot$. We call this regime a self-oscillatory one.

In the case shown in Fig. \ref{1_el_c}, the global attractor in the system is the stable equilibrium. In this case the system can generate a ''response`` to a sufficiently large external impulse. The equations describing the perturbed element can be written in the form
\begin{equation} \label{1_element_ext}
\begin{cases}
\epsilon \mathop{x}\limits^\cdot = x - x^3/3 - y + I(t)\\
\mathop{y}\limits^\cdot = x - a
\end{cases}
\end{equation}
In the simplest case the external stimulus to the element can be modelled using the function $I(t)$ shown in Fig. \ref{1_el_d}. Being perturbed, the state of the element performs some evolution, after which it returns to an equilibrium (see Fig. \ref{1_el_d}). This regime is called \textit{excitable} one.

In the present paper we consider the ensemble of two excitable neuron-like elements of form \eqref{1_element} with symmetric excitatory couplings, which, in accordance with the biological principles  \cite{Destexhe994}, are given by the following function:
\begin{equation} \label{coupling}
I(\phi) = \frac{g}{1 + e^{k(\alpha - \phi)} + e^{k(\phi - \beta)}},
\end{equation}
where $\phi=\arctan\frac{y}{x}$, the parameter $g$ characterizes the strength of the coupling between the elements. For sufficiently large values of the parameter $k$, the coupling function $I(\phi)$, is a smooth and approximates well the function shown in Fig. \ref{1_el_d}. In this formula, $y$ can be approximated by $x$ as follows: if the phase point is in the region of slow motions, then $y\approx x-x^3/3$, otherwise $y=\pm\frac{2}{3}$. 

\begin{figure}[h!]
	\centering
	\subfloat[]
	{
		\includegraphics[width = .4\linewidth]{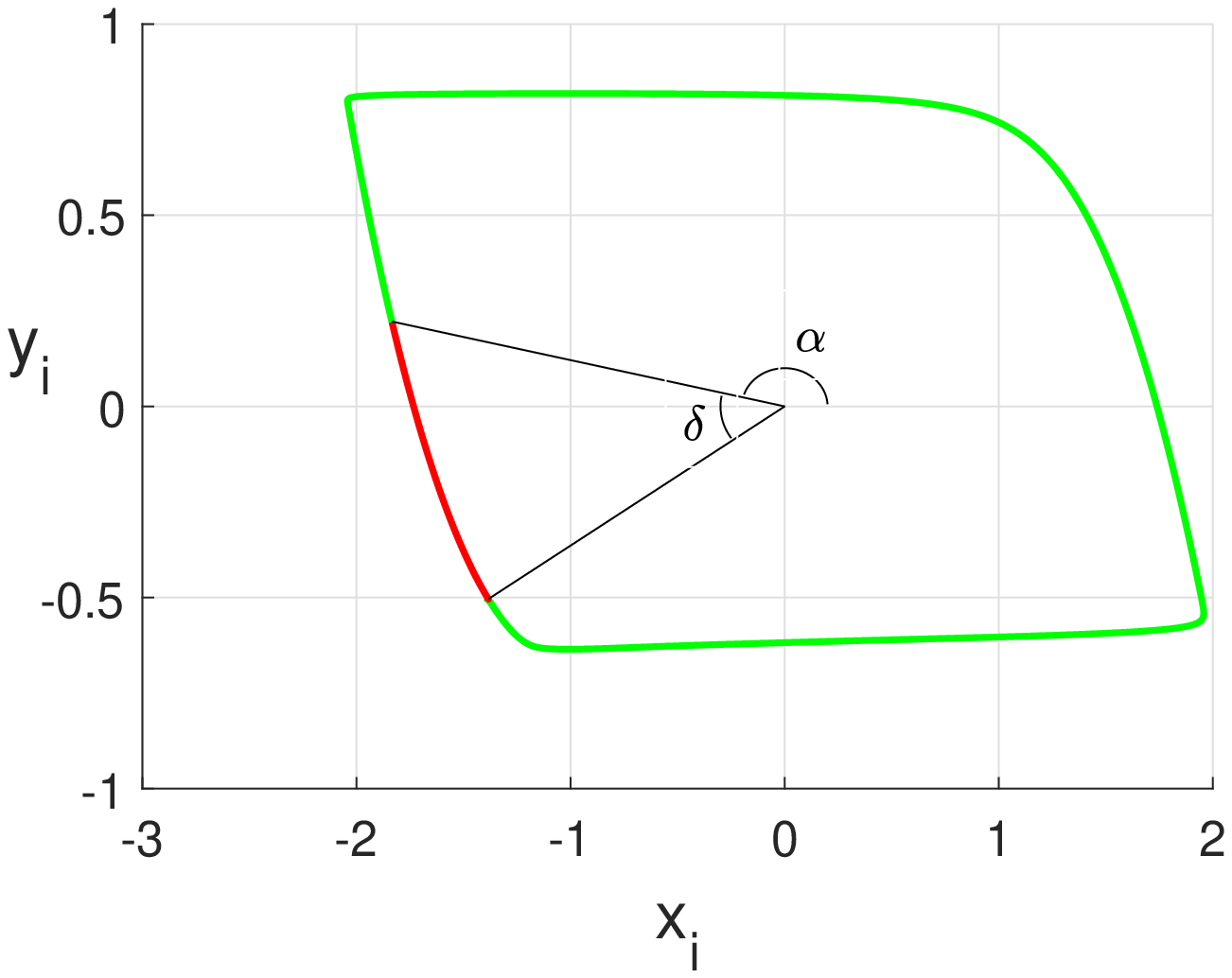}
		\label{mechanism_of_activation_a}
	}
	\subfloat[]
	{
		\includegraphics[width = .4\linewidth]{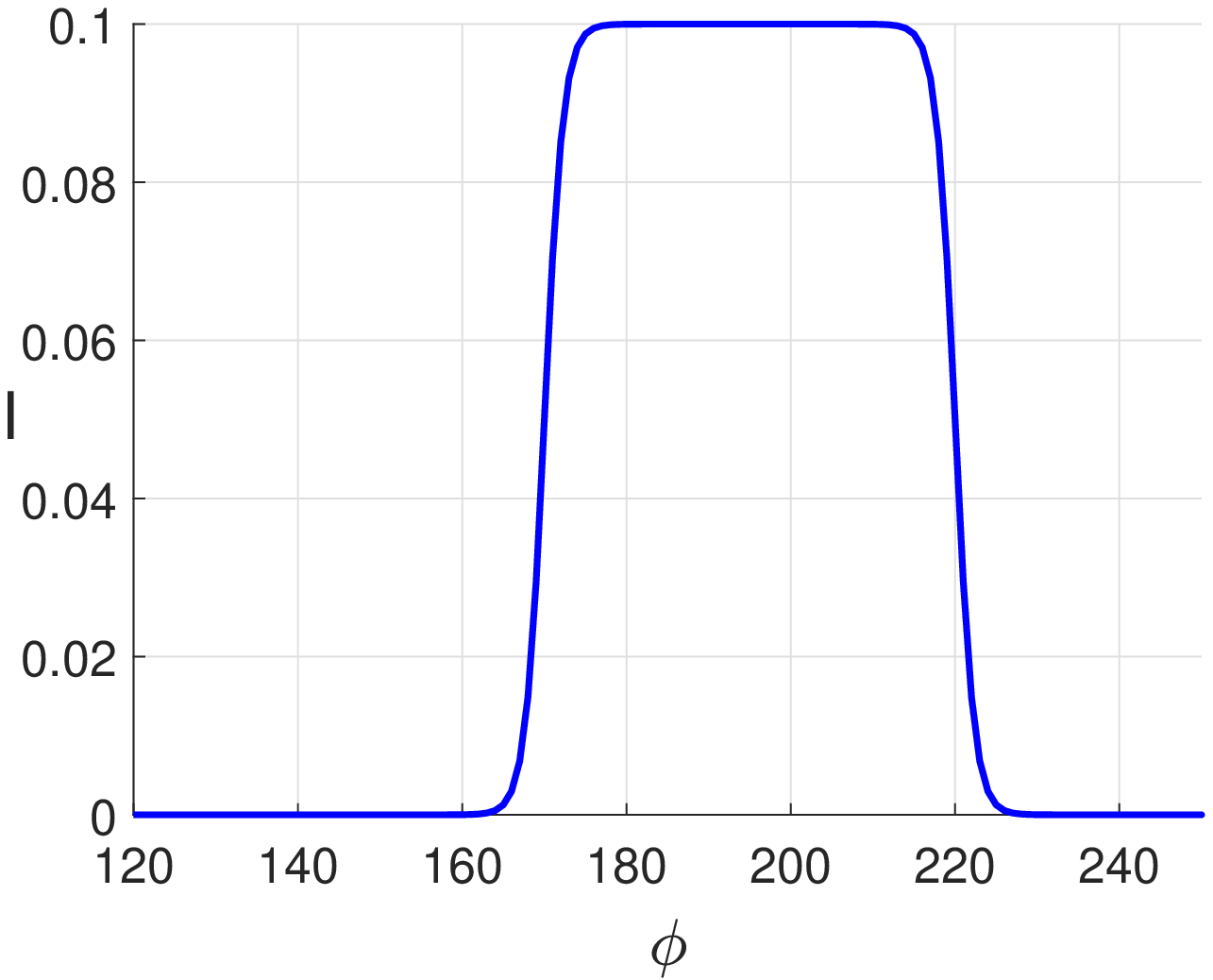}
		\label{mechanism_of_activation_b}
	}
	\caption{(a) The mechanism of transferring activity from one element to another. The projection of the limit cycle onto the phase plane of one element is shown. When the phase point is in the sector between the angles $\alpha$ and $\alpha+\delta$ (part of the projection of the limit cycle intersecting this sector is shown in red), an activating impulse is sent to the other element. (b) Graph of dependence of the activation function $I$ on the phase angle $\phi$.}
	\label{mechanism_of_activation}
\end{figure}

The transmission of activity from one element to another occurs as follows. When the phase $\phi$ of the active presynaptic element reaches $\alpha$, the current of constant amplitude is applied to the postsynaptic element. The time of the impact of this stimulus is defined by the difference $\delta = \beta - \alpha$, which means that the effect is terminated as soon as the representing point of the presynaptic element in the phase plane leaves the sector enclosed between the angles $\alpha$ and $\beta$. If, at the time of activation, the postsynaptic element is in a state close to the resting state, then it will respond. The described mechanism is schematically shown in the Fig. \ref{mechanism_of_activation_a}; the function $I(\phi)$ is shown in Fig. \ref{mechanism_of_activation_b}.

Thus, an ensemble of two coupled neuron-like elements is given by the following system of differential equations:
\begin{equation} \label{ensembles_4eq}
\begin{cases}
\epsilon \mathop{x_1}\limits^\cdot = x_1 - {x_1}^3/3 - y_1 + I(\phi_2)\\
\mathop{y_1}\limits^\cdot = x_1 - a\\
\epsilon \mathop{x_2}\limits^\cdot = x_2 - {x_2}^3/3 - y_2 + I(\phi_1)\\
\mathop{y_2}\limits^\cdot = x_2 - a
\end{cases},
\end{equation}
where $\phi_i = \arctan\frac{y_i}{x_i}$ ($i = 1,2$).

Note that this system is invariant under the change $x_1\leftrightarrow x_2,y_1 \leftrightarrow y_2$. As a result of this symmetry, for each trajectory $(x_2^*(t),y_2^*(t),x_1^*(t),y_1^*(t))$ of the system, there exists the trajectory, symmetric to it with the respect to the invariant plane $\{P:x_1=x_2,y_1=y_2\}$, or this trajectory is self-symmetric (in particular, it lies in the invariant plane $P$).

Further, we fix the following values of the parameters: $a=-1.01$ (the elements are in the excitable regime), $\epsilon=0.01$, $k=50$, $g=0.1$. According to the physical meaning of the parameter $\delta$, which determines the duration of elements activation, $\delta>0$ should be positive and, besides, $\alpha<\beta$. In the next sections we study the impact of the coupling parameters $\alpha$ and $\delta$ on the dynamics of ensemble \eqref{ensembles_4eq}.

\section{Regimes of neuron-like activity}
\begin{figure}[h!]
	\centering
	\subfloat[]
	{
		\includegraphics[width = 0.4\linewidth]{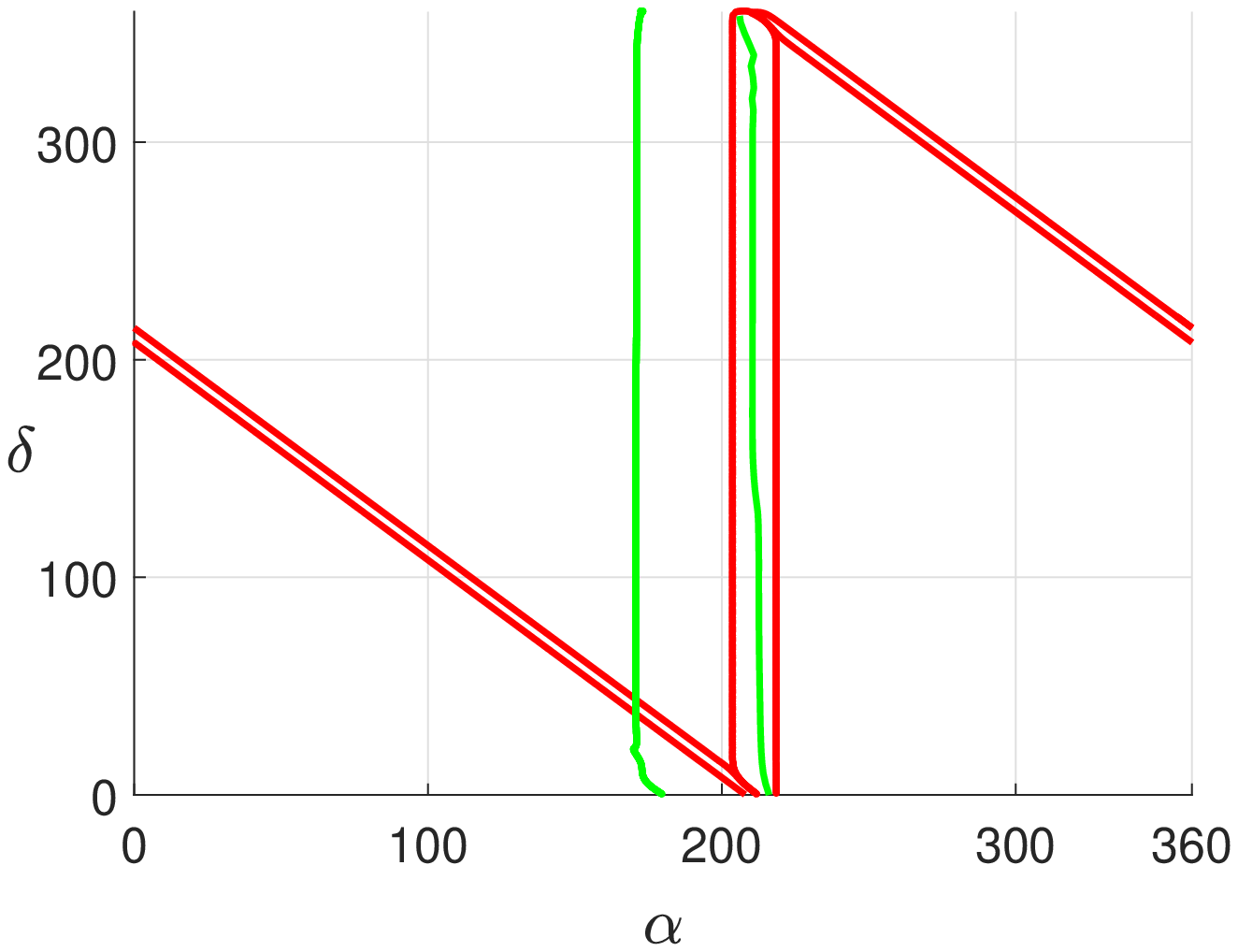}
		\label{bif_diagram_a}
	}
	\subfloat[]
	{
		\includegraphics[width = 0.4\linewidth]{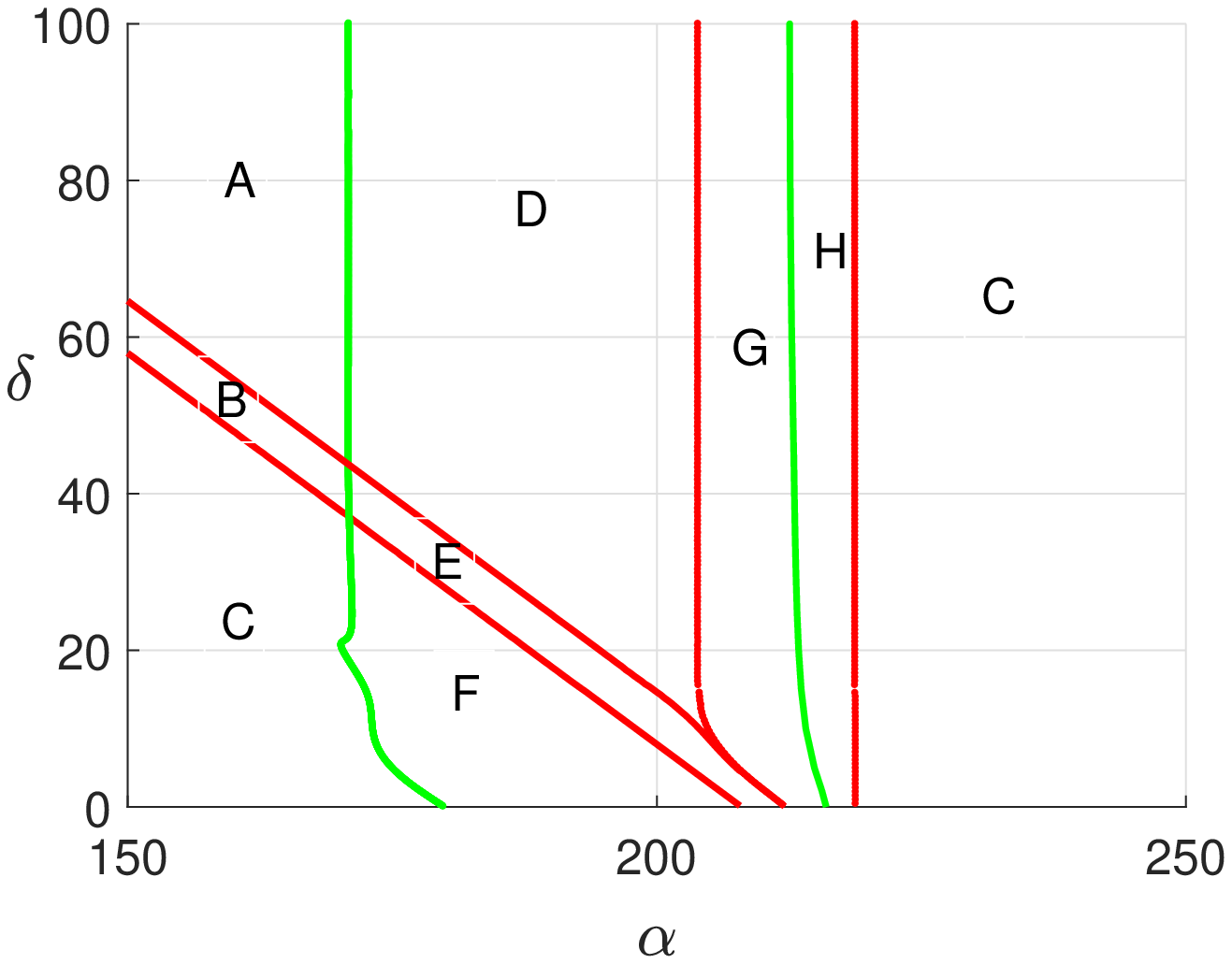}
		\label{bif_diagram_b}
	}
	\caption{Bifurcation diagram and its enlarged fragment. Red colour marks the lines of the Andronov-Hopf bifurcations for the equilibrium, green colour --- lines  of pitchfork bifurcations for limit cycles. In the regions $A$ and $C$ the attractor is a stable equilibrium. In the region $B$ various stable limit cycles corresponding to sequential spiking activity coexist. In the regions $D$ and $F$ stable equilibria and anti-phase limit cycles coexist. In the region $E$ the attractor is a stable anti-phase cycle. In the region $G$ two stable limit cycles coexist: in-phase and anti-phase ones. In the region $H$, one can observe only in-phase stable limit cycle as attractor.}
	\label{bif_diagram}
\end{figure}

\begin{figure}[h!]
	\centering
	\subfloat[$\alpha = 210^\circ$, $\delta = 50^\circ$]
	{
		\includegraphics[width = .3\linewidth]{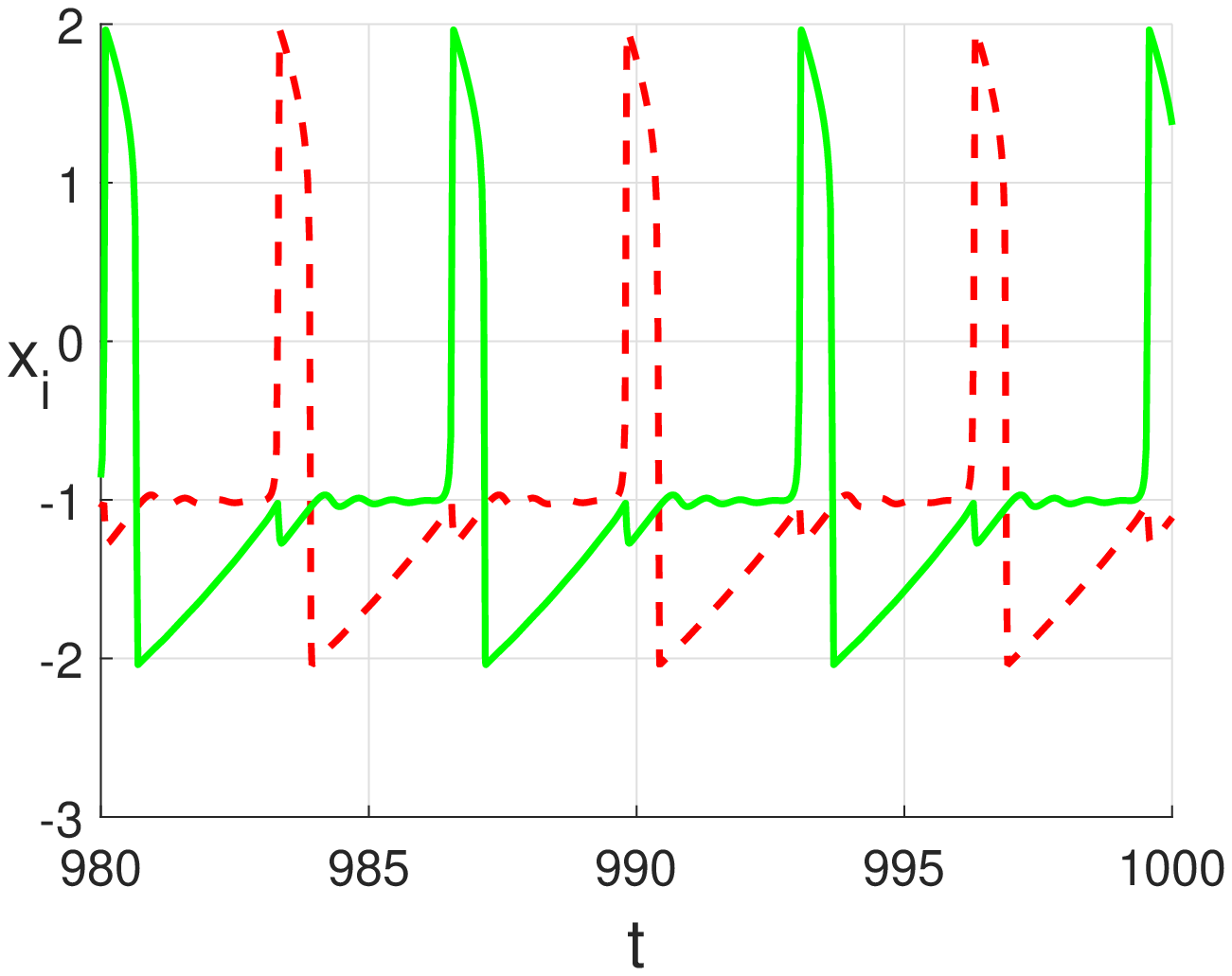}
		\label{time_serieses_a}
	}
	\subfloat[$\alpha = 210^\circ$, $\delta = 50^\circ$]
	{
		\includegraphics[width = .3\linewidth]{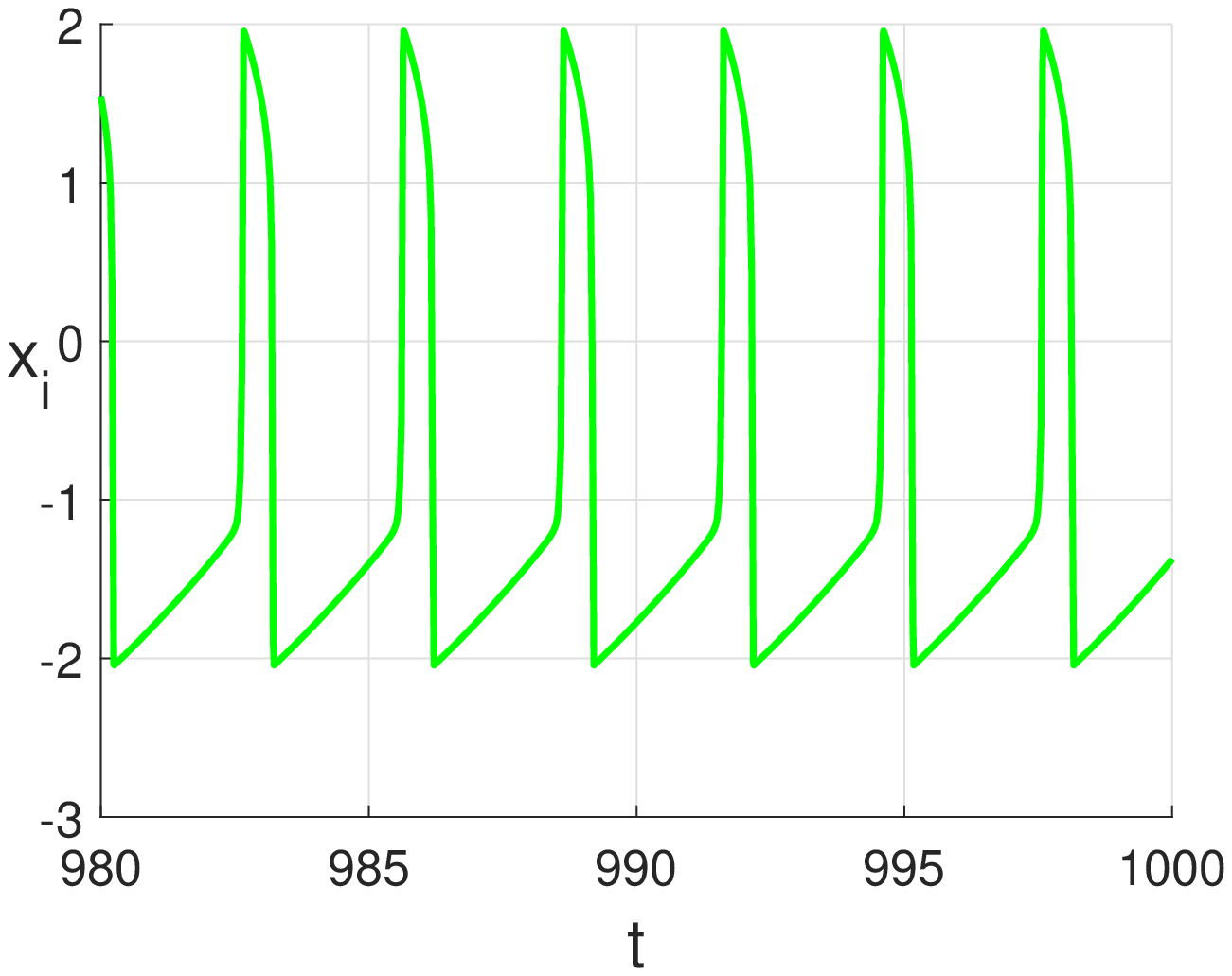}
		\label{time_serieses_b}
	}
	\subfloat[$\alpha = 157^\circ$, $\delta = 50^\circ$]
	{
		\includegraphics[width = .3\linewidth]{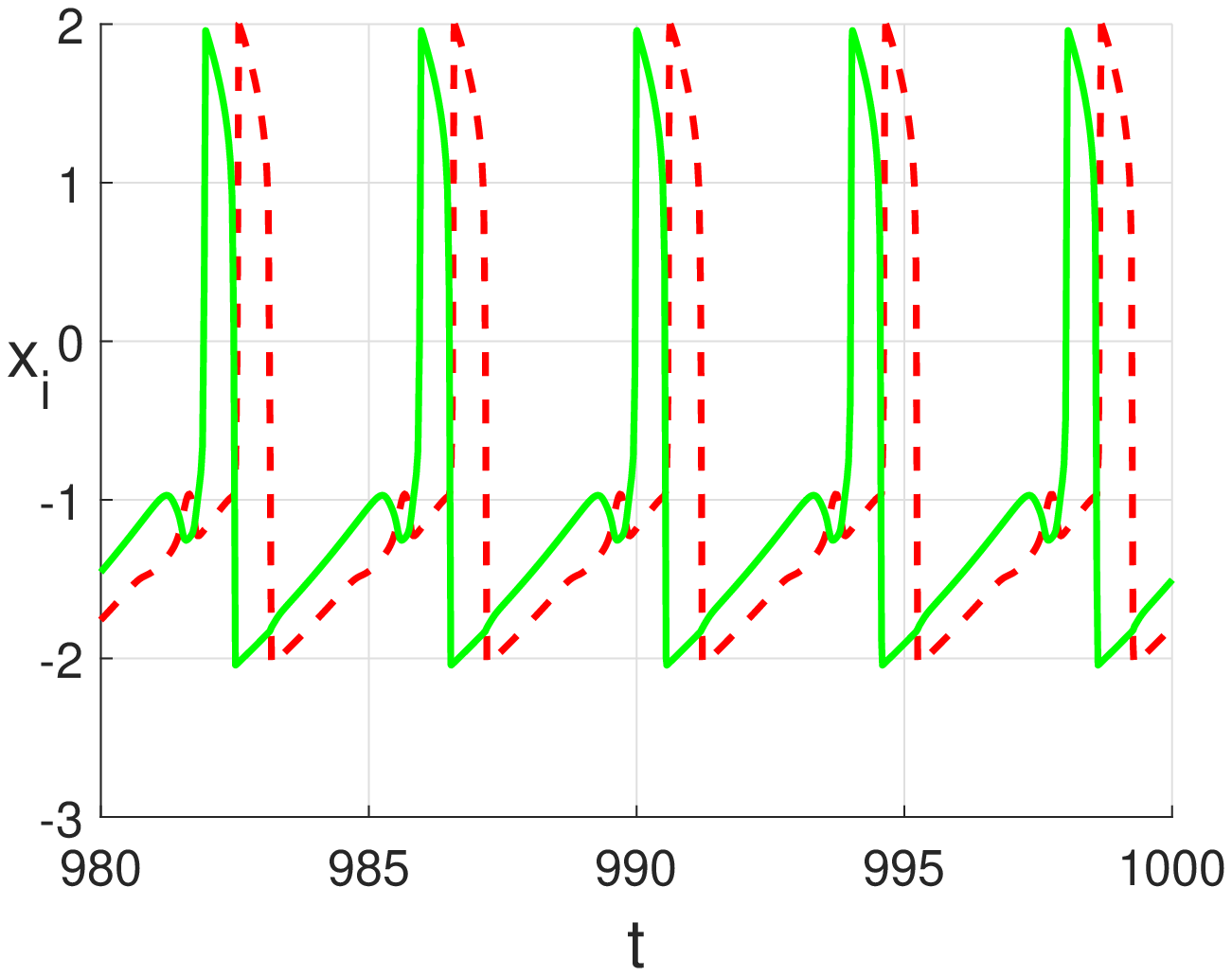}
		\label{time_serieses_c}
	}\\
	\subfloat[$\alpha = 157^\circ$, $\delta = 50^\circ$]
	{
		\includegraphics[width = .3\linewidth]{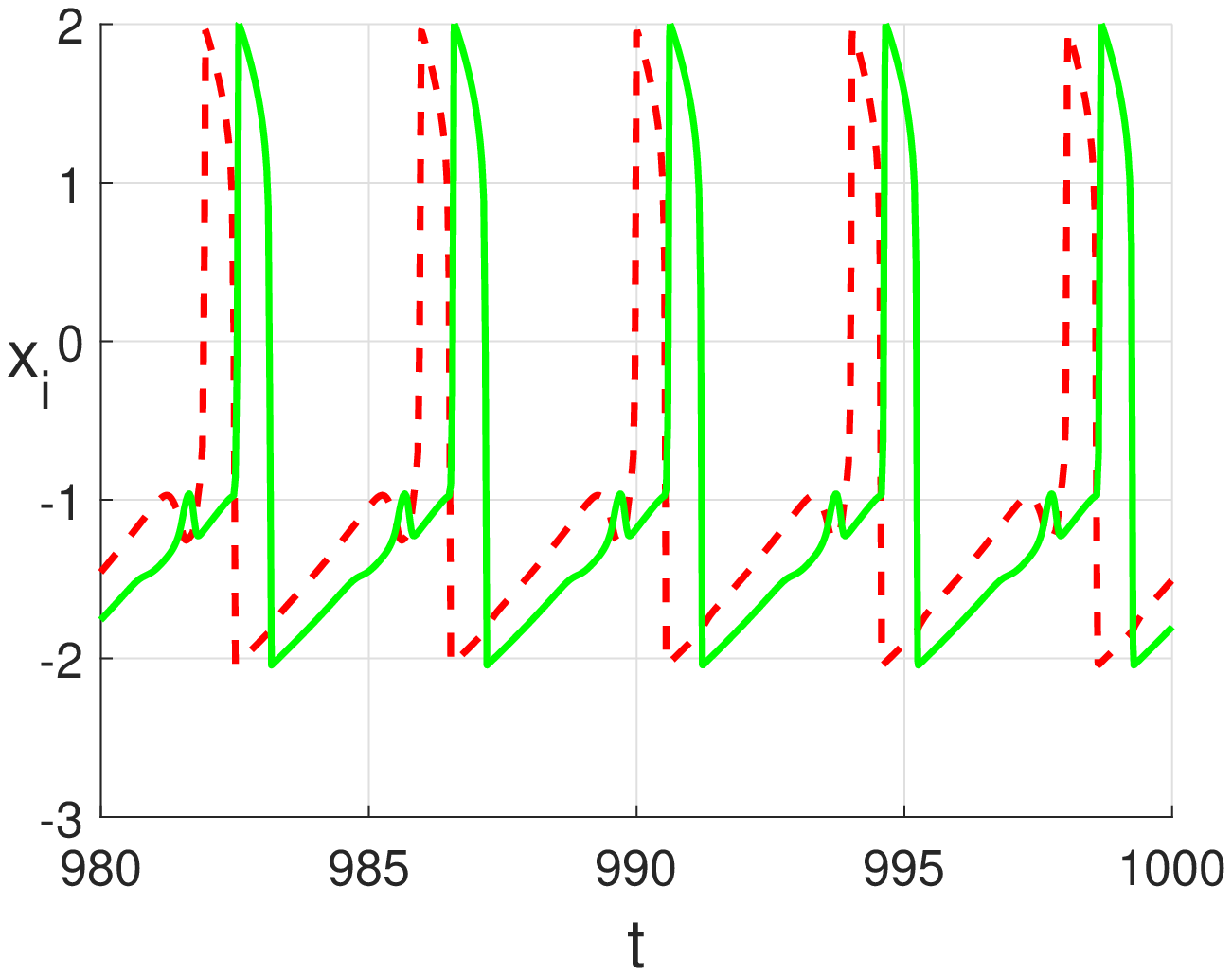}
		\label{time_serieses_d}
	}
	\subfloat[$\alpha = 158^\circ$, $\delta = 50^\circ$]
	{
		\includegraphics[width = .3\linewidth]{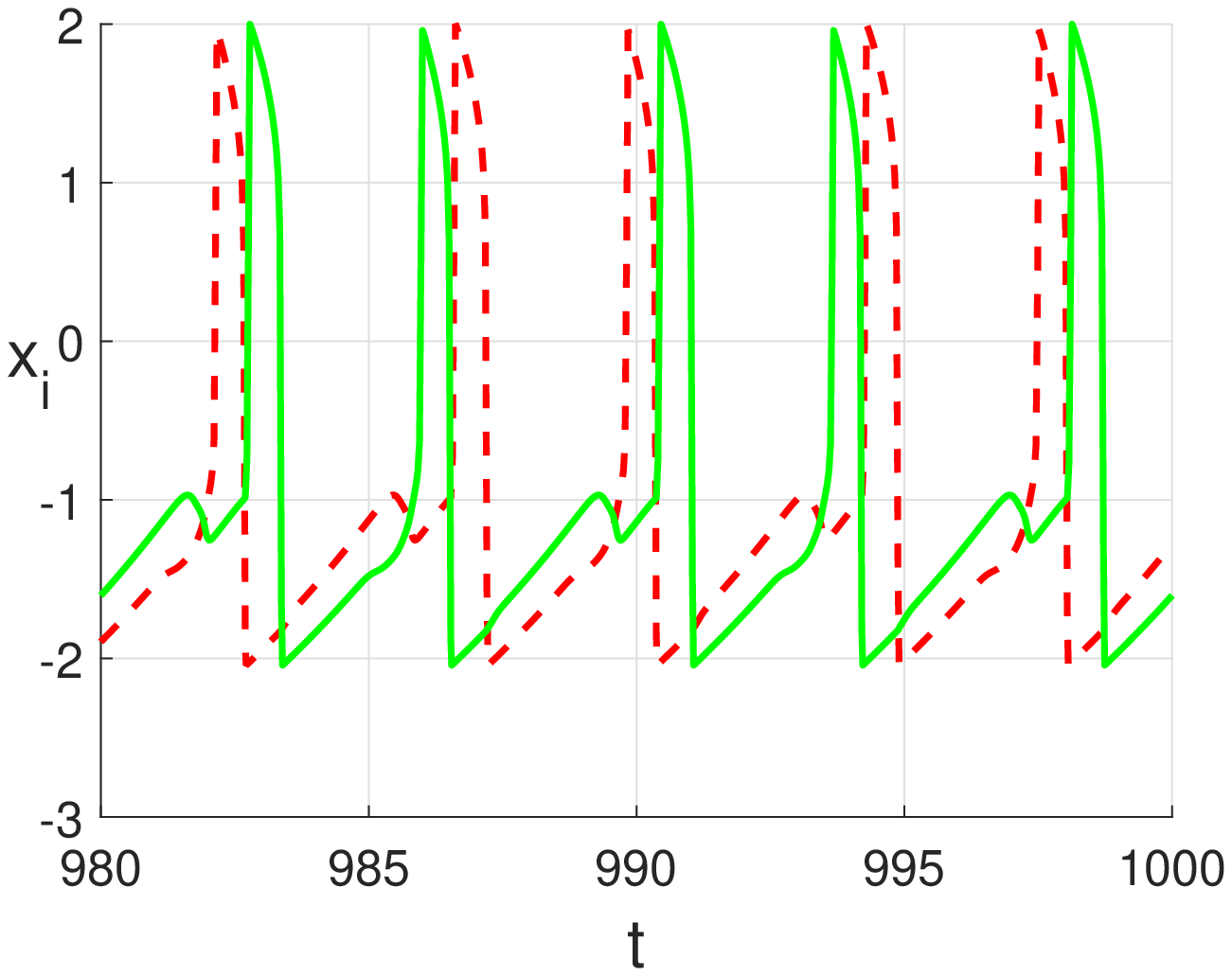}
		\label{time_serieses_e}
	}
	\subfloat[$\alpha = 164.5915^\circ$, $\delta = 50^\circ$]
	{
		\includegraphics[width = .3\linewidth]{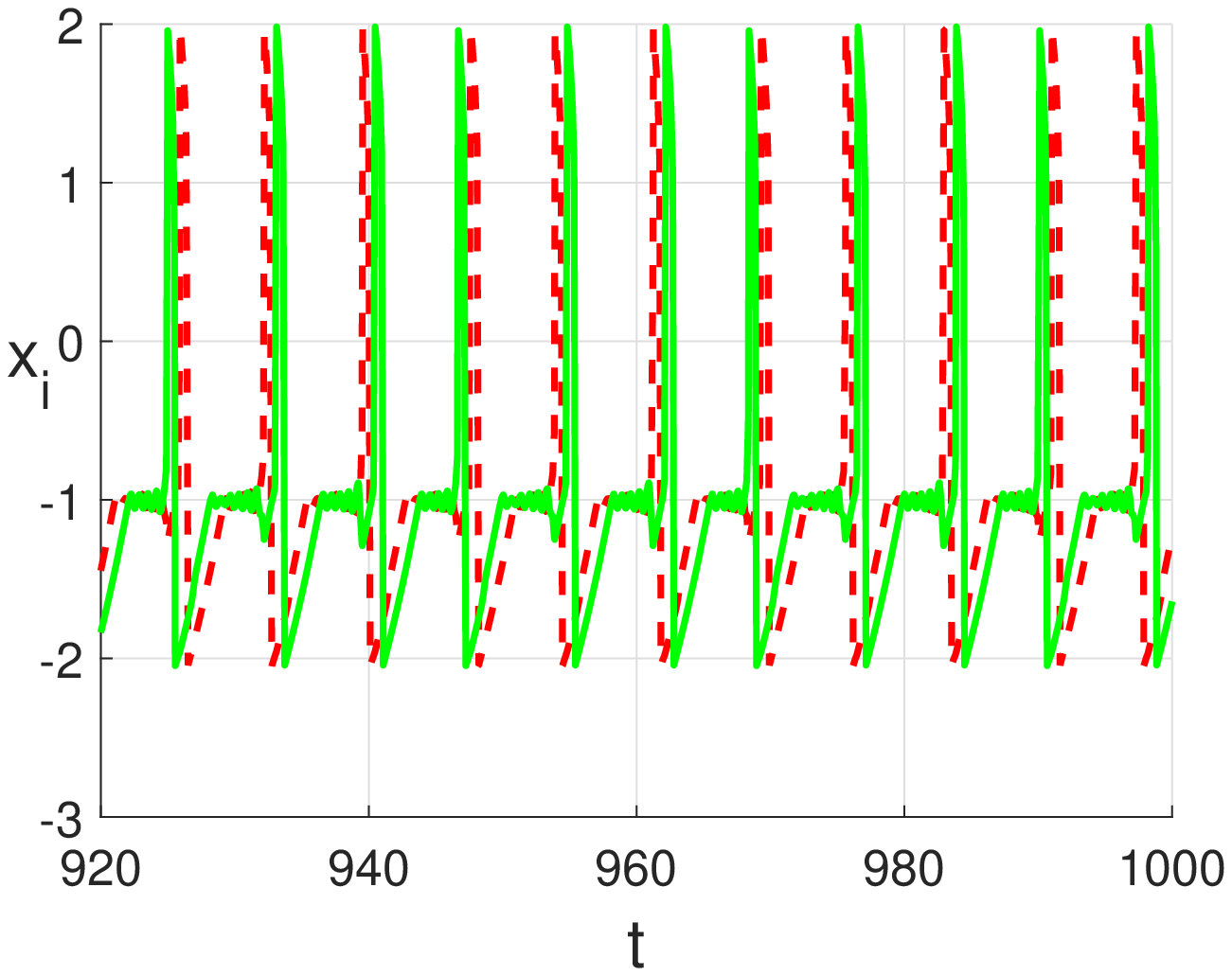}
		\label{time_serieses_f}
	}
	\caption{Time series of the variables $x_1(t)$ and $x_2(t)$ for different regimes. (a) Anti-phase spiking regime $L_{anti}$. (b) In-phase spiking regime $L_{in}$. (c) Sequential spiking activity with the first leading element $L_{12}$. (d) Sequential spiking activity with the second leading element $L_{21}$. (e) Sequential activity with switching order of activation of the elements $L_{1221}$. (f) Sequential activity with complex switchings of the order of activation of elements $L_ {122121}$.}
	\label{time_serieses}
\end{figure}

\begin{figure}[h!]
	\centering
	\subfloat[]
	{
		\includegraphics[width = 0.3\linewidth]{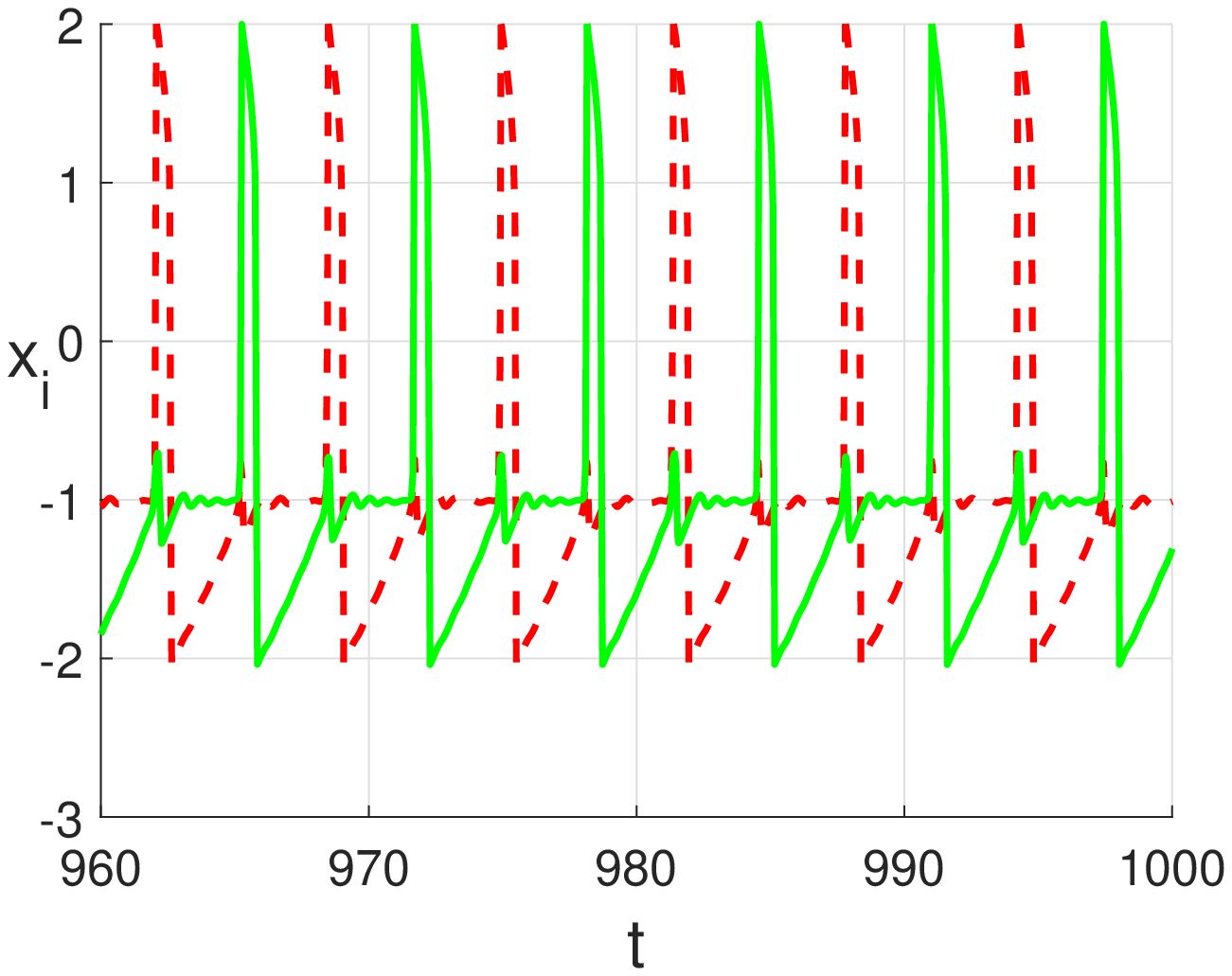}
		\label{time_serieses_chaos_a}
	}
	\subfloat[]
	{
		\includegraphics[width = 0.3\linewidth]{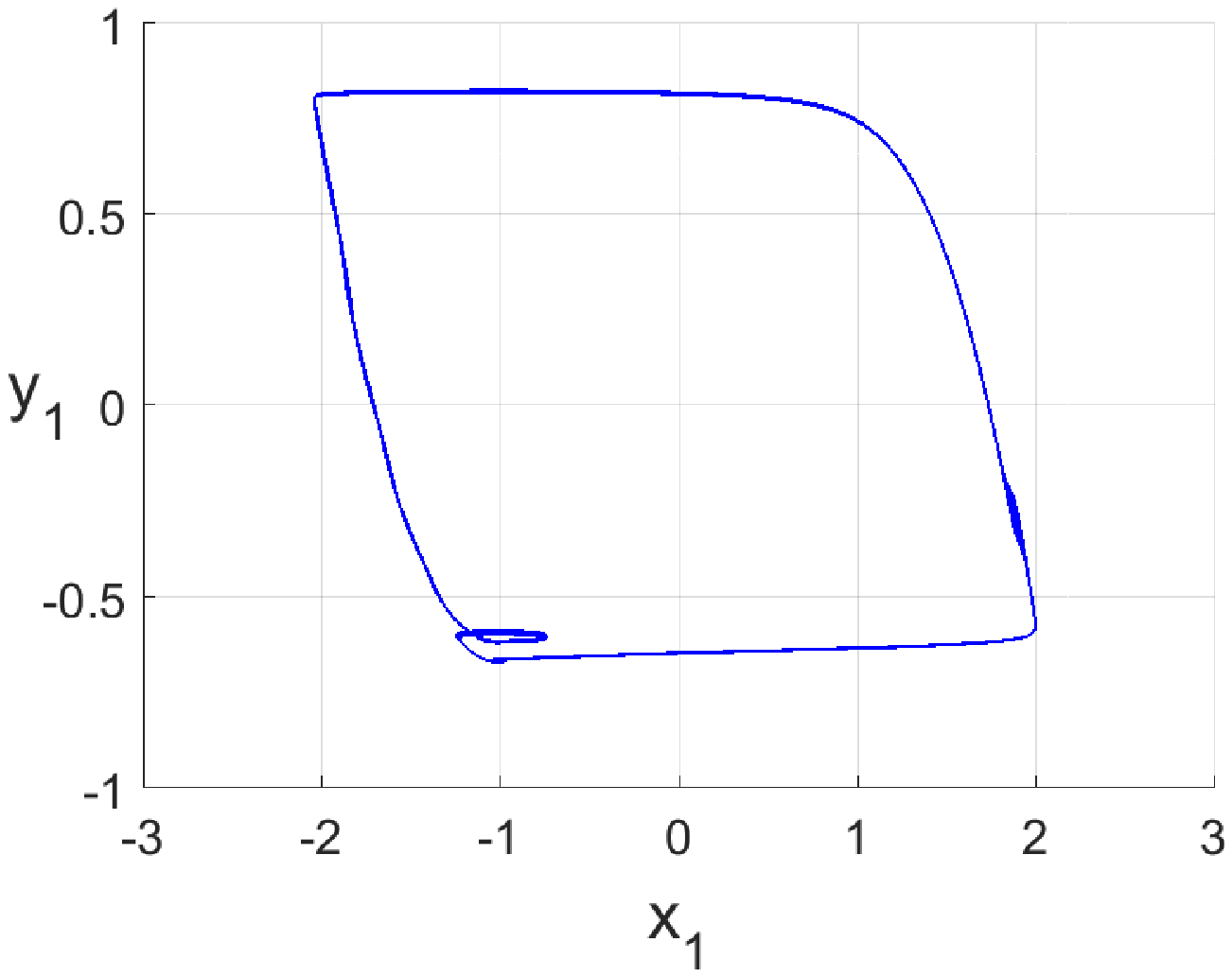}
		\label{time_serieses_chaos_b}
	}
	\subfloat[]
	{
		\includegraphics[width = 0.3\linewidth]{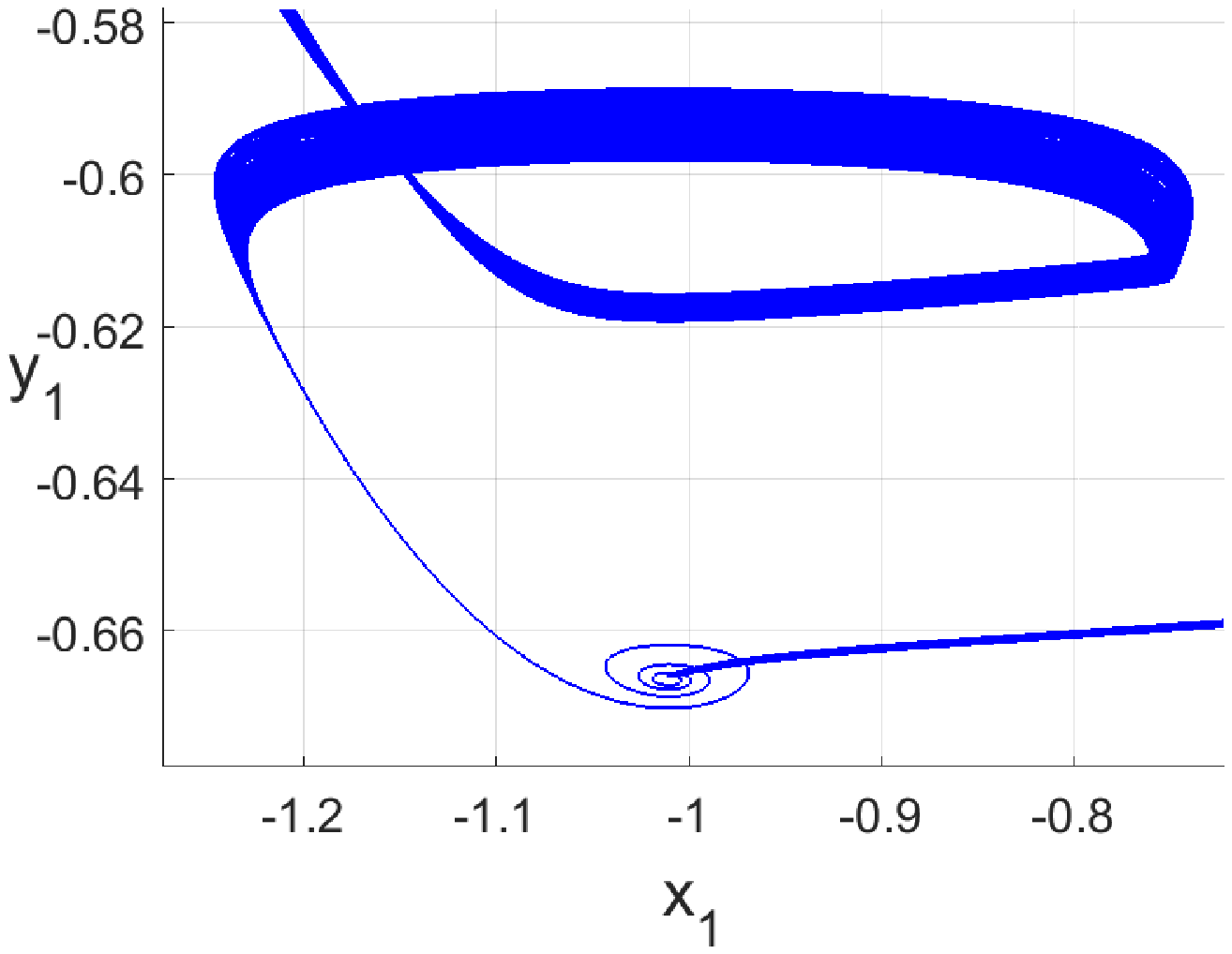}
		\label{time_serieses_chaos_c}
	}
	\caption{(a) Time series and (b), (c) projections of 4-dimension phase space onto the plane $(x_1, y_1)$ for the case of chaotic spiking activity. Here values of the coupling parameters $\alpha = 213.648^\circ$, $\delta = 15^\circ$.}
	\label{time_serieses_chaos}
\end{figure}

It this paper we found that system \eqref{ensembles_4eq}, depending on the control parameters of coupling $\alpha$ and $\delta$, can generate various regimes of neuron-like activity, such as in-phase and anti-phase regimes, various types of sequential spiking activity, including chaotic ones. Using the analytical and numerical methods, that will be described below, the bifurcation diagram shown in Fig. \ref{bif_diagram} was constructed, and the regions of existence of all above described regimes were found on the $(\alpha, \delta)$ parameter plane. In this section, we describe regimes of neuron-like activity that can be observed in each region of the bifurcation diagram. In the next section we will indicate the bifurcations that lead to transitions between described regimes.

First of all, we note that several bistability regions were found in the system. In the regions $D$ and $F$, two regimes coexist: the anti-phase regime (its time series are shown in Fig. \ref{time_serieses_a}) and the regime of quiescence, corresponding to the relatively slow drift towards a stable equilibrium. In the region $G$, there are also two coexisting regimes: the anti-phase regime and the in-phase regime (see the time series in Fig. \ref{time_serieses_b}). The most complex dynamics of the system is observed in the region $B$, where several stable periodic regimes corresponding to sequential activity exist. In this region, we discovered and investigated the sequential activity regimes $L_{12}$ with the leading first element (when the spike in the first element activates the second element, see Fig. \ref{time_serieses_c}) and $L_{21}$ with the leading second element (see Fig. \ref{time_serieses_d}). In addition, there are more complex types of sequential activity in the region $B$, for example, the regimes $L_{1221}$ (see the Fig. \ref{time_serieses_e}) and $L_{122121}$ (see Fig. \ref{time_serieses_f}). Note that we described not all possible types of sequential activity in the region $B$. In forthcoming papers, we plan to investigate in more details bifurcations in this area. All other regions of the bifurcation diagram are characterized by the existence of a single stable limit regime. In regions $A$ and $C$ oscillations (periodic motions) are absent. In this case the attractor is a stable equilibrium. In the regions $E$ and $H$ only existing limit regimes are anti-phase and in-phase ones, respectively.
	
\begin{figure}[h!]
	\centerline{\includegraphics[width = 0.4\columnwidth]{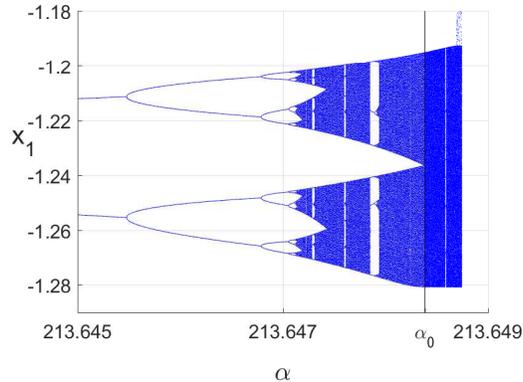}}
	\caption{The coordinate $x_1$ of a fixed point on the Poincare section $y_2=0$ for $\delta=15^\circ$. A strange attractor (the Feigenbaum scenario) appears as a result of the cascade of period-doubling bifurcations.}
	\label{period_doublings}
\end{figure}

In addition to the periodic regimes corresponding to regular activity of various types, a chaotic regime was observed in the system. It appears as a result of the cascade of period doubling bifurcations with the anti-phase limit cycle (see the Fig. \ref{period_doublings}). The time series of the indicated chaotic regime strongly resembles a time series corresponding to the anti-phase limit cycle (compare Fig. \ref{time_serieses_chaos_a} and Fig. \ref{time_serieses_a}). Thus, we call the described chaotic regime the \textit{regime of chaotic anti-phase activity}.
	
\section{Description of the bifurcation diagram}

Here we will give more detailed analysis of the bifurcation diagram of Fig. \ref{bif_diagram} and describe main bifurcations occurred while transitions between regions.

On Fig. \ref{bif_diagram} Andronov-Hopf bifurcations of the equilibrium state correspond to red curves, and the pitch-fork bifurcations of limit-cycles correspond to green ones. The lines $\alpha = 0$ and $\alpha = 360^\circ$ should be identified i. e. the parameter space is a cylinder.

In the regions $A$ and $C$ attractor is the stable equilibrium. The subcritical pitchfork bifurcation takes place while transitions between regions $A$ and $D$ or $C$ and $F$. As result, the saddle anti-phase cycle becomes the stable one. Thus, in regions $D$ and $F$, in addition to the stable equilibrium, another attractor also appears, which is the stable anti-phase cycle $L_{anti}$. The time series of this cycle is shown in Fig. \ref{time_serieses_a}. The transition from region $D$ to region $G$ is accompanied by the supercritical Andronov-Hopf bifurcation (a detailed analysis of this bifurcation can be found in the next section). As result, the stable equilibrium becomes the unstable saddle-focus and the stable in-phase limit cycle $L_{in}$ is born (its time series are shown in Fig. \ref{time_serieses_b}). Thus, in region $G$ two stable limit cycles coexist: in-phase and anti-phase ones. When moving from the region $G$ to the region $H$ the stable anti-phase cycle $L_{anti}$ undergoes supercritical pitchfork bifurcation and it becomes the saddle one and two stable limit cycles $L_{anti1}$ and $L_{anti2}$ are born. Note that the region of coexistence of the limit cycles $L_{anti1}$ and $L_{anti2}$ with the in-phase limit cycle $L_{in}$ is very narrow. Slightly to the right of curve that separate regions $H$ and $G$ the cycles $L_{anti1}$ and $L_{anti2}$ cycles undergo a cascade of period-doubling bifurcations (see Fig. \ref{period_doublings}). As result, two Feigenbaum attractors arise (see region $\alpha < \alpha_0$ in Fig. \ref{period_doublings}) and then these chaotic attractors merge and form single attractor (see region $\alpha > \alpha_0$ in Fig. \ref{period_doublings}). However, the region of existence of this attractor is also very narrow. With a slight increasing parameter $\alpha$ the attractor undergoes crisis after which most trajectories go to the stable in-phase cycle $L_{in}$.
When moving from the region $H$ to the region $C$, the inverse supercritical Andronov-Hopf bifurcation takes place under which the stable limit cycle $L_{in}$ merges with saddle-focus equilibrium which becomes stable one. When moving from regions $D$ or $F$ to the region $E$ the subcritical Andronov-Hopf bifurcation takes place, i. e. the stable equilibrium merges with the saddle limit cycle and becomes saddle-focus one. Thus, in the region $E$ the single attractor (the anti-phase cycle $L_{anti}$) exists. Transition from regions $A$ or $C$ to region $B$ is also accompanied by the subcritical Andronov-Hopf bifurcation.

As it was said in the previous section, the most complicated dynamics is observed in $B$ region. This region can contain several  coexisting stable limit cycles that correspond to various types of sequential activity.

Thus, our researches let us to identify four regions of bistability in space of parameters of the model: $G$ region of coexistence of stable limit cycles $L_{in}$ and $L_{anti}$, region $B$ of coexistence of stable limit cycles $L_{12}$ and $L_{21}$ that correspond to two types of sequential activity, as well as regions $D$ and $F$ of coexistence of the stable equilibrium and the stable limit cycle $L_{anti}$.

\section{Analytical results}
\label{Analytical_results_sect}

Here we present some analytical results on bifurcations of the equilibrium state in system \eqref{ensembles_4eq} and show that for certain values of the parameters an in-phase limit cycle can appear from it as a result of the supercritical Andronov-Hopf bifurcation. The equilibria in the system are determined from the following relations
\begin{equation*}
\begin{cases}
x_1 - {x_1}^3/3 - y_1 + I(\phi_2) = 0\\
x_1 - a = 0\\
x_2 - {x_2}^3/3 - y_2 + I(\phi_1) = 0\\
x_2 - a = 0
\end{cases}
\end{equation*}
Hence, $x_1 = x_2 = a$ and equations for $y_1$ and $y_2$ are as follows:
\begin{equation} \label{equilibrium_2eq}
\begin{cases}
y_1 = a - a^3/3 + I(\arctan \frac{y_2}{a}) = \tilde{I}(y_2)\\
y_2 = a - a^3/3 + I(\arctan \frac{y_1}{a}) = \tilde{I}(y_1)
\end{cases}
\end{equation}
The solutions of this system are fixed points and points of period 2 of the map $\bar{y}=\tilde{I}(y)$. Each fixed point $y_0$ corresponds to the equilibrium state $O (a, y_0, a, y_0)$ of the system \eqref{ensembles_4eq}, and a pair of points of period 2 $y_{10}=\tilde{I}(y_{20})$ and $y_{20}=\tilde{I}(y_{10})$ corresponds to the the pair of equilibria $O_1(a, y_ {10}, a, y_ {20})$ and $O_2(a, y_ {20}, a, y_ {10})$.
The Jacobi matrix of the system \eqref{ensembles_4eq} has the following form
\begin{equation*}
\begin{bmatrix}
\frac{1 - {x_1}^2}{\epsilon} & -\frac{1}{\epsilon} & \frac{\partial I(\phi_2)}{\epsilon \partial x_2} & \frac{\partial I(\phi_2)}{\epsilon \partial y_2}\\
1 & 0 & 0 & 0 \\
\frac{\partial I(\phi_1)}{\epsilon \partial x_1} & \frac{\partial I(\phi_1)}{\epsilon \partial y_1} & \frac{1 - {x_2}^2}{\epsilon} & -\frac{1}{\epsilon}\\
0 & 0 & 1 & 0
\end{bmatrix}
\end{equation*}
Consequently, the characteristic equation can be written in the form
\begin{equation} \label{characteristic_equation}
\left(\lambda \left(\lambda - \frac{1 - {x_1}^2}{\epsilon}\right) + \frac{1}{\epsilon}\right)\left(\lambda \left(\lambda - \frac{1 - {x_2}^2}{\epsilon}\right) + \frac{1}{\epsilon}\right) - \left(\frac{I_{x_1}}{\epsilon} \lambda + \frac{I_{y_1}}{\epsilon}\right) \left(\frac{I_{x_2}}{\epsilon} \lambda + \frac{I_{y_2}}{\epsilon}\right) = 0
\end{equation}
In the equilibrium $O(a,\,y_0,\,a,\,y_0)$ for partial derivatives of the coupling function $I(\phi)$, the following relations hold: $I_{x_1}(a,y_0)=I_{x_2}(a,y_0)=I_x$, $I_{y_1}(a,y_0)=I_{y_2}(a,y_0)=I_y$. Then the characteristic equation \eqref{characteristic_equation} can be written in the form
\begin{equation*}
\left(\lambda \left(\lambda - \frac{1 - a^2}{\epsilon}\right) + \frac{1}{\epsilon}\right)^2 - \left(\frac{I_x}{\epsilon} \lambda + \frac{I_y}{\epsilon}\right)^2 = 0.
\end{equation*}
Let us find roots of this characteristic equation:
\begin{equation*}
\begin{cases}
\lambda_{1,2} = \frac{1 - a^2 + I_x \pm \sqrt{(1 - a^2)^2 + 2(1 - a^2)I_x + {I_x}^2 - 4\epsilon(1 - I_y)}}{2\epsilon}\\
\lambda_{3,4} = \frac{1 - a^2 - I_x \pm \sqrt{(1 - a^2)^2 - 2(1 - a^2)I_x + {I_x}^2 - 4\epsilon(1 - I_y)}}{2\epsilon}
\end{cases}
\end{equation*}
Thus, the equilibrium $O(a,\,y_0,\,a,\,y_0)$ undergoes Andronov-Hopf bifurcation when one of the two following conditions is satisfied:
\begin{equation} \label{condition1}
\begin{cases}
1 - a^2 + I_x = 0\\
(1 - a^2)^2 + 2(1 - a^2)I_x + {I_x}^2 - 4\epsilon(1 - I_y) < 0
\end{cases}
\end{equation}
or
\begin{equation} \label{condition2}
\begin{cases}
1 - a^2 - I_x = 0\\
(1 - a^2)^2 - 2(1 - a^2)I_x + {I_x}^2 - 4\epsilon(1 + I_y) < 0
\end{cases}
\end{equation}
Conditions \eqref{condition1} and \eqref{condition2} can be written in the following form:
\begin{equation}
\label{AH_cond_1}
\left \lbrace 
\begin{array}{c}
1 - a^2 + I_x = 0\\
I_y < 1 
\end{array}
\right . \text{or}\;
\left \lbrace
\begin{array}{c}
1 - a^2 - I_x = 0\\
I_y > -1
\end{array}
\right .
\end{equation}
Now we define the partial derivatives $I_x$ and $I_y$ of the coupling function given by the equation \eqref{coupling}:
\begin{equation} \label{deriv}
\begin{cases}
\frac{\partial I}{\partial x} = \frac{gky\left(e^{k(\arctan \frac{y}{x} - \beta)} - e^{k(\alpha - \arctan \frac{y}{x})}\right)}{\left(x^2 + y^2\right)\left(1 + e^{k(\alpha-\arctan \frac{y}{x})} + e^{k(\arctan \frac{y}{x} - \beta)}\right)^2}\\
\frac{\partial I}{\partial y} = -\frac{gkx\left(e^{k(\arctan \frac{y}{x} - \beta)} - e^{k(\alpha - \arctan \frac{y}{x})}\right)}{\left(x^2 + y^2\right)\left(1 + e^{k(\alpha-\arctan \frac{y}{x})} + e^{k(\arctan \frac{y}{x} - \beta)}\right)^2} = -\frac{x}{y}\frac{\partial I}{\partial x}
\end{cases}
\end{equation}
Using relation \eqref{deriv} and conditions \eqref{AH_cond_1}, one can obtain the following expressions:
\begin{equation}
\label{AH_conds}
\left \lbrace 
\begin{array}{c}
I_x = a^2 - 1  \\
\frac{a(1 - a^2)}{y_0} < 1  
\end{array}
\right . \text{or}\;
\left \lbrace
\begin{array}{c}
I_x = 1 - a^2 \\
\frac{a(1 - a^2)}{y_0} < 1 
\end{array}
\right .
\end{equation}
Since $0<I(\phi)<g$, it follows that $a-\frac{a^3}{3}<\tilde{I}(y)<a-\frac{a^3}{3}+g$. Then it follows from the system \eqref{equilibrium_2eq} that $a-\frac{a^3}{3}<y_0<a-\frac{a^3}{3}+g$. With the chosen values of the parameters $a$ and $g$, it follows from the last inequality that $y_0<0$, and hence the inequality $\frac{a(1-a^2)}{y_0}<1$ is fulfilled. As a result, in the conditions of \eqref{AH_conds} only equalities should be considered. Adding these equalities to the relation $a-a^3/3-y_0+I(\arctan \frac{y_0}{a})=0$ for determining the coordinate $y_0$ of the equilibrium state $O$, we obtain two conditions that define the Andronov-Hopf bifurcations:
\begin{equation} \label{AH_cond_1_end}
\begin{cases}
I_x = a^2 - 1\\
a - a^3/3 - y_0 + I(\arctan \frac{y_0}{a}) = 0
\end{cases}
\end{equation}
or
\begin{equation} \label{AH_cond_2_end}
\begin{cases}
I_x = 1 - a^2\\
a - a^3/3 - y_0 + I(\arctan \frac{y_0}{a}) = 0
\end{cases}
\end{equation}

Next, we determine the conditions for the appearance of an in-phase cycle by qualitative methods. For the in-phase limit cycle, the conditions $x_1(t)=x_2(t)=x(t)$ and $y_1(t)=y_2(t)=y(t)$ are satisfied. Then the system \eqref{ensembles_4eq} can be written in the form
\begin{equation*}
\begin{cases}
\epsilon \mathop{x}\limits^\cdot = x - x^3/3 - y + I(\phi)\\
\mathop{y}\limits^\cdot = x - a
\end{cases}
\end{equation*}
For sufficiently small $I(\phi)$, the implicit function $y(x)$ given by
\begin{equation} \label{imlicit_func}
F(x, y) = x - x^3/3 - y + I(\phi) = 0,
\end{equation}
has 2 extrema. The bifurcation resulting in the generation of the in-phase cycle occurs when the function $y(x)$ has an extremum at $x=a$: $y'(a)=0$ (the birth of the limit cycle occurs when the minimum point moves to the left from the line $x=a$, otherwise the cycle disappears) (see Fig. \ref{1_el}). Let us find the derivative of the function implicitly defined by the relation \eqref{imlicit_func}: $\frac{dy}{dx}=\frac{1-x^2+I_x}{1-I_y}$. As a result, the bifurcation curve of the cycle birth is found from condition
\begin{equation*}
\begin{cases}
1 - a^2 + I_x = 0\\
a - a^3/3 - y_0 + I(\arctan \frac{y_0}{a}) = 0
\end{cases},
\end{equation*}
which coincides with the condition \eqref{AH_cond_1_end} of the Andronov-Hopf bifurcation.

In order to obtain the condition \eqref{AH_cond_2_end} from the same considerations, one should note that for coupling function $I(\phi; \alpha, \beta)$ the following equality is satisfied: $I(\phi; \alpha, \beta) \approx g - I(\phi; \beta, \alpha + 2\pi)$ (equality $\lim\limits_{k \to +\infty}\left(I(\phi; \alpha, \beta) + I(\phi; \beta, \alpha + 2\pi)\right) = g$ also is satisfied). Then the implicit function $y(x)$ is given by the equality $x-x^3/3-y+g-I(\phi)=0$, and its derivative has the form $\frac{dy}{dx} = \frac{1 - x^2 - I_x}{1 + I_y}$. Then the bifurcation curve of the cycle birth will be given by the system
\begin{equation*}
\begin{cases}
1 - a^2 - I_x = 0\\
a - a^3/3 - y_0 + I(\arctan \frac{y_0}{a}) = 0
\end{cases}
\end{equation*}
This condition is the same as the condition \eqref{AH_cond_2_end}.

Thus, we shown that, an in-phase limit cycle appears as a result of the Andronov-Hopf bifurcation, and bifurcation curves are analytically found.

\section{Conclusions}
The phenomenological model of the ensemble of two excitable FitzHugh-Nagumo elements with symmetric excitatory couplings is proposed. The novelty of the proposed model is that here a new type of coupling constructing by means of  a smooth function approximating a standard rectangular function was proposed.

Despite its simplicity, the proposed model demonstrates a rich variety of different regimes of neuron-like activity. So, we observed regular regimes of in-phase, anti-phase and sequential activity. Each of these regimes corresponds to a stable periodic motion of a certain type in the phase space of the system. In addition to regimes of regular activity, the system can also demonstrate the anti-phase regime of chaotic activity corresponding to the strange attractor that appears as a result of the cascade of period doubling bifurcations.

Moreover, the detailed studies of bifurcations leading to the appearance of these regimes have been carried out. On the plane of the parameters specifying the beginning $\alpha$ and duration $\delta$ of the couplings impact, the regions corresponding to different types of neuron-like activity are determined. It is shown that there are regions of bistability for which the system has different stable regimes. For example, in some parameter regions, depending on the choice of the initial conditions, the ensemble can demonstrate either the anti-phase or in-phase periodic regimes. There are also regions of parameters where different types of sequential activity coexist.

In the future, we plan to apply the developed phenomenological model for modeling neural ensembles from a large number of elements.

\begin{acknowledgments}
We thank Prof. S. V. Gonchenko for valuable advices and comments. The study was supported by RSF grant 17-72-10228 (sections 1-3) and RFBR grant 17-02-00467 (section 4). A.O. Kazakov acknowledges the support of the Program of  fundamental studies of NRU HSE in 2018.
\end{acknowledgments}

\def\cprime{$'$}


\begin{thebibliography}{34}%

\bibitem{Markram}
{Underwood E. Science. 2013. V~339. P.~1022.}

\bibitem{DAngelo}
{D'Angelo E., Solinas S., Garrido J., et al. Functional Neurology. 2013. V.~28, №~3. P. 153.}

\bibitem{Synapses}
{Cowan W.M., Sudhof T.C., Stevens C.F, editors. Synapses. Baltimore: The Johns Hopkins University Press, 2001.}

\bibitem{Neuroscience}
{Purves D., Augustine G.J., Fitzpatrick D., et al., editors. Neuroscience. 2nd edition. Sunderland (MA): Sinauer Associates, 2001.}

\bibitem{Models1}
{Izhikevich E. M. Dynamical systems in neuroscience. MIT press, 2007.}

\bibitem{Details}
{Koch C. Biophysics of Computation: Information Processing in Single Neurons. New York: Oxford University Press, 1999;
	De Schutter E., editor. Computational Neuroscience: Realistic Modeling for Experimentalists. Boca Raton, FL: CRC Press, 2000. }

\bibitem{HH}
{Hodgkin A.L., Huxley A.F. The Journal of Physiology. 1952. V.~4, N.~117. P.~500.}

\bibitem{Shilnikov}
{Malaschenko T., Shilnikov A., Cymbalyuk G. Physics Review E. 2011. V. 84. P. 041910.}

\bibitem{Shilnikov1}
{Malaschenko T.*, Shilnikov A., Cymbalyuk G. PLoS ONE. 2011. V. 6, N.7. P. e21782.}

\bibitem{Shilnikov2}
{Cymbalyuk G., Calabrese R., Shilnikov A. Proc. of the Annual Computational Neuroscience Meeting. 2003}

\bibitem{FHN}
{Schwan H.P., editor. Biological Engineering. New York: McGraw-Hill Book Co., 1969. P.~1-85).}

\bibitem{Roth-Rossum}
{De Schutter E., editor. Computational Modeling Methods for Neuroscientists. Cambridge: MIT Press, 2009. 432 p.}

\bibitem{Baladron-Fasoli-Faugeras-Touboul}
{Baladron J., Fasoli D., Faugeras O., Touboul J. Journal of Mathematical Neuroscience. 2012. V.~2. P.~10.}

\bibitem{Gao}
{Keener J., Sneyd J. Mathematical Physiology. New York: Springer, 1998.}

\bibitem{Mobius}
{Valles-Codina O., Mobius R., Rudiger S., Schimansky-Geier L. Physical Review E. 2011. V.~83. P.~036209.}

\bibitem{Rankovic}
{Rankovic D. Matematicki Vesnik. 2011. V.~63, N~2. P.103.}

\bibitem{Shin-Lee-Kim}
{Shin CW., Lee SG.,  Kim S. Journal of the Korean Physical Society. 2006. V.~48. P.~179.}

\bibitem{Hansel-Mato-Meuner}
{Hansel D., Mato G., Meunier C. Europhysics Letters. 1993. V.~23, N.~5. P.~367.}

\bibitem{Izhikevich2000}
{Izhikevich E.M., Int. J. Bifurcation and Chaos. 2000. V.~10, N.6. P. 1171.}

\bibitem{Doss-Bachelet-Francoise-Piquet}
{Doss-Bachelet C., Francoise JP., Piquet C. ComPlexUs. 2003. V.~1, N.~3. P.~101.}

\bibitem{Komarov-Osipov-Suykens-Rabinovich}
{Komarov M.A., Osipov G.V., Suykens J.A.K., Rabinovich M.I. CHAOS. 2009. V.~19. P.~015107.}

\bibitem{Rabinovich}
{Rabinovich M., Volkovskii A., Lecanda P., Huerta R., Abarbanel H. D. I., Laurent G. Physical Review Letters. 2001. V.~87, N.~6. P.~068102.}

\bibitem{Nguyen}
{Nguyen L.H., Hong K.S. Mathematics and Computers in Simulation. 2011. V.~82. P.~590.}

\bibitem{Binczak}
{Binczak S., Jacquir S., Bilbault J.M., Kasantsev V., Nekorkin V. Neural Networks. 2006. V.~19, N.~5. P.~684.}

\bibitem{Jacquir}
{Jacquir S., Binczak S., Bilbault J.M., Kazantsev V., Nekorkin V. Nonlinear Dynamics. 2006. V.~44. P.~29.}

\bibitem{Hoff}
{Hoff A., dos Santos J.V., Manchein C., Albuquerque H.A. The European Physical Journal B. 2014. V.~87, N.~7. P.~151.}

\bibitem{Campbell}
{Campbell S.A., Waite M. Nonlinear Analysis. 2001. V.~47. P.~1093.}

\bibitem{Tehrani}
{Tehrani N.F., Razvan M. Mathematical Biosciences. 2015. V.~270. P.~41.}

\bibitem{Yanagita}
{Yanagita T., Ichinomiya T., Oyama Y. Physical Review E. 2005. V.~72. P.~056218.}

\bibitem{Brown-Feng-Feerick}
{Brown D., Feng J., Feerick S. Physical Review Letters. 1999. V.~82, N.~23. P.~4731.}

\bibitem{Wang-Zhang-Deng}
{Wang J., Zhang T., Deng B. Chaos, Solitons and Fractals. 2007. V.~31. P.~30.}

\bibitem{Song}
{Song Y., Xu J. IEEE Transactions on Neural Networks and Learning Systems. 2012. V.~23, N.~10. P.~1659.}

\bibitem{Wang-Tian-Dhamala-Liu}
{Wang Z., Tian C., Dhamala M., Liu Z. Scientific Reports. 2017. V.~7. P.~561.}

\bibitem{Destexhe994}
{Destexhe A., Meinen Z.F., Sejnowski T.J. Neur. Comput. 1994. V.~6. P.~14.}
\end{thebibliography}
\end{document}